\documentclass[12pt]{amsart}%
\usepackage{amsfonts}
\usepackage{txfonts}
\usepackage{mathtools}
\usepackage{amsmath}
\usepackage{amssymb}
\usepackage{amsthm}
\usepackage{newlfont}
\usepackage{graphicx}%
\providecommand{\U}[1]{\protect\rule{.1in}{.1in}}
\newtheorem{thm}{Theorem}[section]
\newtheorem{theorem}{Theorem}[section]

\newtheorem{cor}[thm]{Corollary}
\newtheorem{lem}[thm]{Lemma}
\newtheorem{prop}[thm]{Proposition}

\theoremstyle{definition}

\theoremstyle{remark}

\numberwithin{equation}{section}

\newcommand{\ed}{\end {document}}

\begin{document}
\title[Hamiltonian stationary manifolds]{On the regularity of Hamiltonian stationary Lagrangian manifolds}
\author{Jingyi Chen}
\address{Department of Mathematics\\
The University of British Columbia.}

\email{jychen@math.ubc.ca}
\author{Micah Warren}
\address{
Department of Mathematics\\
University of Oregon, Eugene, OR 97403, U.S.A.}
\email{micahw@uoregon.edu}
\thanks{The first author was supported in by NSERC Discovery Grant 22R80062. The
second author was partially supported by NSF Grant DMS-1438359. }

\begin{abstract}
We prove a Morrey-type theorem for Hamiltonian stationary Lagrangian submanifolds of
$\mathbb{C}^{n}$: If a $C^{1}$ Lagrangian submanifold is a critical point of
the volume functional under Hamiltonian variations, then it must be real
analytic. Locally, a Hamiltonian stationary manifold is determined
geometrically by harmonicity of its Lagrangian phase function, or
variationally by a nonlinear fourth order elliptic equation of the potential
function whose gradient graph defines the Hamiltonian stationary submanifolds
locally. Our result shows that Morrey's theorem for minimal submanifolds
admits a complete fourth order analogue. We establish full regularity and
removability of singular sets of capacity zero for weak solutions to the
fourth order equation with $C^{1,1}$ norm below a dimensional constant, and to
$C^{1,1}$ potential functions, under certain convexity conditions, whose
Lagrangian phase functions are weakly harmonic.

\end{abstract}
\maketitle

\section{Introduction}

In this paper, we study regularity of Hamiltonian stationary submanifolds of
complex Euclidean space. These are critical points of the volume functional
under Hamiltonian variations, and locally they are governed by a fourth order
nonlinear elliptic equation. We show, among other results, that when a
Hamiltonian stationary manifold is $C^{1}$, then it must be real analytic. For
minimal submanifolds, a classical theorem of Morrey states: If a minimal
submanifold of Euclidean space is $C^{1},$ then it is real analytic
\cite[Theorem 10.7.1]{MorreyBook}. Our approach to the fourth order equation
is completely different from Morrey's for the second order minimal surface
equations. Our result applies when the fourth order equation is satisfied away
from a set of capacity zero. This echoes the extendibility results of
\cite[Theorem 1.2]{HL75}, where it is shown that solutions to the system of
minimal surface equations on a domain in $\mathbb{R}^{n}$ extend across closed
sets of zero $(n-1)$-dimensional Hausdorff measure.

We now describe the analytic setup of the geometric variational problem. For a
fixed bounded domain $\Omega\subset{\mathbb{R}}^{n}$, let $u:\Omega
\rightarrow\mathbb{R}$ be a smooth function. The gradient graph $\Gamma
_{u}=\left\{  \left(  x,Du(x)\right)  :x\in\Omega\right\}  $ is a Lagrangian
$n$-dimensional submanifold in $\mathbb{C}^{n}$, with respect to the complex
structure $J$ defined by the complex coordinates $z_{j}=x_{j}+\sqrt{-1}y_{j}$
for $j=1,\cdots,n$. The volume of $\Gamma_{u}$ is given by
\[
F_{\Omega}(u)=\int_{\Omega}\sqrt{\det\left(  I+(D^{2}u)^{T}D^{2}u\right)
}dx.
\]
A twice differentiable function $u$ is critical for $F_{\Omega}(u)$ under
compactly supported variations of the scalar function $u$ if and only if $u$
satisfies the Euler-Lagrange equation%
\begin{equation}
\int_{\Omega}\sqrt{\det g}g^{ij}\delta^{kl}u_{ik}\eta_{jl}\,dx=0\text{
}\ \ \ \text{for all }\eta\in C_{c}^{\infty}(\Omega). \label{vHamstat}%
\end{equation}
Here, summation convention is applied over repeated indices, $\delta^{kl}$ is
the Kronecker delta, and $g$ is the induced metric from the Euclidean metric
on $\mathbb{R}^{2n}$, which can be written as
\[
g=I+(D^{2}u)^{T}D^{2}u.
\]
We can define the volume $F_{\Omega}(u)$ whenever $u\in W^{2,n}(\Omega),$ so
$W^{2,n}(\Omega)$ is a natural space on which to seek critical points. We will
call \eqref{vHamstat} {\it the variational Hamiltonian stationary equation}. A
function $u\in W^{2,n}(\Omega)$ is called a weak solution the variational
Hamiltonian equation if \eqref{vHamstat} holds.

If the potential $u$ is in $C^{4}(\Omega)$, the equation (\ref{vHamstat}) is
equivalent to the following \textit{geometric Hamiltonian stationary
equation}
\begin{equation}
\Delta_{g}\theta=0 \label{Hamstat}%
\end{equation}
where $\Delta_{g}$ is the Laplace-Beltrami operator on $\Gamma_{u}$ for the
induced metric $g$ (cf. \cite{MR1202805}, \cite[Proposition 2.2]{SW}). The
function $\theta$ is called the Lagrangian phase function for the gradient
graph $\Gamma_{u}$ and is defined by
\[
\theta=\operatorname{Im}\log\det\left(  I+\sqrt{-1}D^{2}u\right)
\]
or equivalently,
\begin{equation}
\theta=\sum_{i=1}^{n}\arctan\lambda_{i} \label{thetadef2}%
\end{equation}
for $\lambda_{i}$ the eigenvalues of $D^{2}u.$ The mean curvature vector along
$\Gamma_{u}$ can be written
\[
\vec{H}=-J\nabla\theta
\]
where $\nabla$ is the gradient operator of $\Gamma_{u}$ for the metric $g,$
see (\cite[2.19]{HL}).
We say a function $u$ is a weak solution of (\ref{Hamstat}) if

\begin{enumerate}
\item The Hessian $D^{2}u$ is defined almost everywhere and $u\in
W^{2,n}(\Omega).$

\item The quantity $\theta$ in (\ref{thetadef2}) is in $W^{1,2}(\Omega).$

\item For all $\eta\in C_{c}^{\infty}(\Omega)$%
\begin{equation}
\int_{\Gamma_{u}}\langle\nabla\theta,\nabla\eta\rangle d\mu_{g}=0.
\label{weakintegral}%
\end{equation}

\end{enumerate}


From an elliptic PDE point of view, the equation (\ref{Hamstat}) is much
preferred: The equation (\ref{Hamstat}) is a second order operator upon a
second order quantity, so we may use the full power of the well-developed
second order nonlinear elliptic theory against the equation. Importantly, the
function (\ref{thetadef2}) is a concave quantity when $\theta$ falls in
certain ranges, or when $u$ is convex. On the other hand, nonlinear double
divergence equations of the form (\ref{vHamstat}) are not as well understood.
We will compare the geometric settings of the two equations in more depth in
Section 2.

A smooth Lagrangian submanifold $L\subset\mathbb{C}^{n}$ that solves
(\ref{Hamstat}) is called \textit{Hamiltonian stationary}. Note that one can
always define the Lagrangian phase function $\theta$, up to an additive
constantt $2k\pi$. In general, a Hamiltonian stationary submanifold in a
symplectic manifold is a critical point of the volume functional under
Hamiltonian deformations, that is, the variations generated by $J\nabla\eta$
for some smooth compactly supported function $\eta$ on ${\mathbb{C}}^{n}$.
Recall that if $u$ satisfies \textit{the special Lagrangian equation}
\cite{HL}
\begin{equation}
\nabla\theta=0 \label{sLag}%
\end{equation}
i.e. $\vec{H}\equiv0$, then the submanifold is critical for the volume
functional under \textit{all} compactly supported variations of the surface
$\Gamma_{u}$. The special Lagrangians are Hamiltonian stationary. The Clifford
torus in the complex plane is Hamitonian stationary but not special
Lagrangain. There are non-flat cones that are Hamiltonian stationary but not
special Lagrangian, and this regularity issue causes serious problems for
constructing minimal Lagrangian surfaces in a K\"ahler-Einstein surface (see
\cite{SW}).

Hamiltonian stationary submanifolds form an interesting class of Lagrangians
in a symplectic manifold as critical points of the volume functional under
Hamiltonian deformations. They generalize the minimal Lagrangian submanifolds
in a K\"ahler-Einstein manifold, especially, the special Lagrangians in a
Calabi-Yau manifold. The existence and stability problem has been studied by
many people via different approaches (cf \cite{MR1062973}, \cite{MR1611051},
\cite{SWJDG}, \cite{HeleinRomon2002}, \cite{MR1986315}, \cite{HeleinRomon2005}%
, \cite{JLS}, and references therein). Yet, a general theory for existence
remains open.

Our first goal is to study the regularity of submanifolds that locally are
described by potentials satisfying (\ref{vHamstat}). In particular, we will
show that if $D^{2}u$ does not have large discontinuities then the potential
$u$ must be smooth, hence solving both \eqref{vHamstat} and \eqref{Hamstat}.
We will consider regularity for weak solutions that lie in the Sobolev space
$W_{loc}^{2,\infty}(\Omega)$.

\begin{thm}
\label{T1}Let $\Omega$ be a domain in $\mathbb{R}^{n}$ and let $Q\subset
\Omega$ be a compact subset (possible empty) with capacity zero. There is a
$c(n)>0$ such that if $u\in C^{1,1}(\Omega\backslash Q)$ is a weak
solution to (\ref{vHamstat}) on $\Omega\backslash Q$ satisfying
\[
\left\Vert u\right\Vert _{C^{1,1}(\Omega\backslash Q)}\leq c(n),
\]
then $u$ is a smooth solution of both (\ref{Hamstat}) and (\ref{vHamstat}) on
$\Omega.$
\end{thm}

Recall that the capacity of a set $Q$ is defined as
\[
\mbox{Cap}(Q)=\inf_{\substack{\phi\in C_{c}^{\infty}(%
\mathbb{R}
^{n}),\\0\leq\phi\leq1,\\\phi=1\text{ near }Q}}\int\left\vert D\phi\right\vert
^{2}dx.
\]
In particular, if the Hausdorff dimension of $Q$ is less than $n-2$ then
Cap($Q$) is zero.

We make several remarks: First, by a rotation, one can choose a gradient graph
representation of $\Gamma$ so that $D^{2}u(0)=0,$ at any point where the
tangent space is defined. Next, as\ there are no size restrictions on $\Omega
$, any continuity condition on the Hessian will suffice. More details are
provided in section 3. Finally, this $c(n)$ is not obtained by a compactness
argument, and can be made explicit.

Next we show that in certain cases where a (slightly weaker) Hessian bound is
assumed, weak solutions to (\ref{Hamstat}) enjoy full regularity.

\begin{thm}
\label{T2}Suppose that $u\in C^{1,1}\left(  \mathbb{B}_{1}(0)\right)  $ and
$u$ is a weak solution of (\ref{Hamstat})$.$ \ If either
\begin{equation}
\theta\geq\delta+\frac{\pi}{2}(n-2)\ \ \text{ a.e.} \label{delta1}%
\end{equation}
for some constant $\delta\in(0,\pi)$; or
\begin{equation}
u-\delta\frac{\left\vert x\right\vert ^{2}}{2}\ \text{is convex }
\label{delta2}%
\end{equation}
for some constant $\delta>0$; or
\begin{equation}
\left\Vert u\right\Vert _{C^{1,1}(\mathbb{B}_{1}(0))}\leq1-\delta
\label{delta3}%
\end{equation}
for some constant $\delta\in(0,1)$, then for $k\geq2$ we have
\[
\left\Vert u\right\Vert _{C^{k,\alpha}(\mathbb{B}_{1/2}(0))}\leq
C(k,n,\left\Vert u\right\Vert _{C^{1,1}({\mathbb{B}_{1}}(0))},\delta).
\]

The conclusion still holds if $\mathbb{B}_{1}(0)$ is replaced by
$\mathbb{B}_{1}(0) \backslash Q$, where $Q$ is a compact subset of
$\mathbb{B}_{1}(0)$ with capacity zero.
\end{thm}

Our strategy is as follows: For a weak solution $u$ to equation
\eqref{Hamstat}, if $\Vert u\Vert_{C^{1,1}}$ is strictly below 1, then the
Lewy-Yuan rotation, adapted to the non-smooth setting (see Proposition
\ref{weakrotation}), converts the question to the case that a (new) potential
function is uniformly convex, that is, (\ref{delta2}), and then the machinery
of viscosity solutions for concave operators applies. Note that the situation
(\ref{delta1}) can be dealt with using the same concave operator theory.
Essentially, this is the Schauder theory for concave equations in \cite{CC}
applied to the inhomogeneous equation of special Lagrangian type. For
extending solutions across $Q$, we invoke a removability theorem of Serrin
\cite{Serrin} for equations in divergence form. For a weak solution $u$ to
\eqref{vHamstat} with small $C^{1,1}$ norm, first we show that $u$ is in
$W_{loc}^{3,2}$, and this allows approximations by smooth functions in
$W_{loc}^{3,2}$ norm and then leads to that $\theta$ (which is a priori merely
$L^{\infty})$ satisfies (\ref{weakintegral}), therefore, the full regularity
obtained for equation \eqref{Hamstat} applies.

To prove our main geometric result, we combine the above two theorems as
follows. Choosing an appropriate tangent plane, locally, we apply Theorem
\ref{T1}. Since the equation (\ref{Hamstat}) is geometrically invariant (up to
an immaterial additive constant), we may rotate the coordinates to where the
quantity $\theta$ is concave, and apply Theorem \ref{T2} to obtain a
description of smoothness of the same manifold. We have

\begin{thm}
\label{T3} Any $C^{1}$ Hamiltonian stationary submanifold of ${\mathbb{C}}%
^{n}$ is real analytic. More generally, suppose $u\in$ $W^{2,n}\left(
\Omega\right)  $, and $u$ satisfies equation (\ref{vHamstat}) on $\Omega$.
There is a constant $c_{0}\left(  n\right)  $ such that if the image of the
tangent planes (where defined) of the gradient graph
\[
\Gamma_{u}=\left\{  \left(  x,Du(x)\right)  :x\in\Omega\right\}
\]
lies in a ball of radius $c_{0}(n)$ in the Grassmannian $Gr(n,2n)$, then
$\Gamma_{u}$ is a real analytic submanifold of $\mathbb{R}^{2n}$.
\end{thm}

In particular, if $D^{2}u$ is within distance $c(n)$ to a continuous function,
then $u$ must be smooth, hence real analytic. For example, while we cannot
rule out non-flat tangent cones occurring, we can rule out non-flat tangent
cones that are nearly flat.

In two dimensions, regularity results have been obtained by Schoen and Wolfson
\cite[Theorem 4.7]{SW} in a general K\"{a}hler manifold setting, where
singularities are known to occur. The examples of singularities are
non-graphical over an open domain \cite[Section 7]{SWJDG}. On the other hand,
the Euclidean case of \cite[Proposition 4.6]{SWJDG} states that $u$ solving
(\ref{Hamstat}) is smooth whenever $u\in C^{2,\alpha}$. Our Theorem \ref{T3}
is a generalization of this result, see Corollary \ref{C^2}.


The rest of the paper is organized as follows. In section 2, we derive and
compare the Euler-Lagrange equations, given mild regularity conditions on $u.$
\ In section 3, we show that nonlinear divergence type fourth order equations
enjoy a regularity boost from $W^{2,\infty}$ to $W^{3,2}$ given a condition on
the nonlinearity, and from this prove Theorem \ref{T1}. In section 4, we give
details on the Lewy-Yuan rotation, as this will be necessary to prove the
third part of Theorem \ref{T2}. In section 5, we discuss and apply the
Schauder theory for equations of special Lagrangian type, showing Schauder
type results when the equation is concave. We then prove Theorem \ref{T2}
under the first two conditions and combine this with the results from section
4 to give us the result in the third case. Theorem \ref{T3} will follow.

\section{Derivation of the Euler-Lagrange equations}

Consider the functional on the space of $C^{2}$ functions on a bounded domain
$\Omega$ in $\mathbb{R}^{n}$
\begin{equation}
F_{\Omega}(u)=\int_{\Omega}\sqrt{\det\left(  I+\left(  D^{2}u\right)
^{T}D^{2}u\right)  }dx. \label{V1}%
\end{equation}
Note that for the gradient graph of a function $u$, we have the induced metric%
\begin{equation}
g_{ij}=\delta_{ij}+u_{ik}\delta^{kl}u_{lj} \label{metricgij}%
\end{equation}
in which case the above functional becomes
\begin{equation}
F_{\Omega}(u)=\int_{\Gamma_{u}}\sqrt{\det g}dx. \label{V2}%
\end{equation}

\begin{prop}
\bigskip Suppose that $u\in C^{3}(\Omega)$. Then $u$ is a weak solution to
(\ref{vHamstat}) on $\Omega$ if and only if $u$ is a weak solution to
(\ref{Hamstat}) on $\Omega$, in which case (\ref{vHamstat}) and (\ref{Hamstat}%
) are each the Euler-Lagrange equation for the functional (\ref{V1}).
\end{prop}

\begin{proof}
First we consider the case where $u$ solves (\ref{vHamstat}). Take a variation
generated by $\eta\in C_{c}^{\infty}(\Omega)$, which varies the manifold along
the $y$-direction in $\mathbb{C}^{n}.$ Computing the volume for the path of
potentials%
\begin{equation}
\gamma\lbrack t](x)=u(x)+t\eta(x), \label{verteta}%
\end{equation}
we get%
\begin{align*}
\left.  \frac{d}{dt}F_{\Omega}(\gamma\lbrack t])\right\vert _{t=0}  &
=\int_{\Omega}\frac{1}{2}\sqrt{g[t]}g^{ij}[t]\left.  \frac{d}{dt}%
g_{ij}[t]\right\vert _{t=0}dx\\
&  =\frac{1}{2}\int_{\Omega}\sqrt{g}g^{ij}\left(  u_{ik}\delta^{kl}\eta
_{lj}+\eta_{ik}\delta^{kl}u_{lj}\right)  dx\\
&  =\int_{\Omega}\sqrt{g}g^{ij}u_{ik}\delta^{kl}\eta_{lj}dx.
\end{align*}
Thus, the first variation of $F_{\Omega}$ at $u$ is given by
\[
\delta F_{\Omega}(\eta)=\int_{\Omega}\sqrt{g}g^{ij}u_{ik}\delta^{kl}\eta
_{lj}dx.
\]
We note that while defining $F_{\Omega}(u)$ requires only that $u\in
W^{2,n}(\Omega)$.


On the other hand, we may compute the variation using the standard first
variational formula for (\ref{V2}) , when $u\in C^{3}$:
\[
\left.  \frac{d}{dt}F_{\Omega}(\gamma\lbrack t])\right\vert _{t=0}=\frac
{d}{dt}Vol(\Gamma_{u})=\int_{\Omega}\langle-\vec{H},V\rangle d\mu_{g}%
\]
where $\vec{H}$ is the mean curvature vector, and $V$ is the variational
field. Recall that the variation $V$ is Hamiltonian if $V=JDf$ for some
compactly supported function $f$ in ${\mathbb{C}}^{n}$. For a Lagrangian
submanifold, we also have \cite[2.19]{HL},
\[
\vec{H}=-J\nabla\theta.
\]
Therefore, a $C^{2}$ Lagrangian submanifold is critical for the volume
functional under Hamiltonian variations if and only if its Lagrangian phase is
weakly harmonic.

In our case, namely, the gradient graph of $u\in C^{3}(\Omega)$, we have a
vertical variational field that is Hamiltonian:
\begin{equation}
V(x)=\left.  \frac{d}{dt}(x,Du(x)+tD\eta(x))\right\vert _{t=0}=(0,D\eta(x)).
\label{vfd}%
\end{equation}

We claim that $u$ is a weak solution to \eqref{Hamstat} is equivalent to that
the gradient graph is critical for all vertical variations. In fact,
\begin{align*}
\delta F_{\Omega}(\eta)  &  =\int_{\Omega}\langle J\nabla\theta,(0,D\eta
)\rangle d\mu_{g}\\
&  =\int_{\Omega}\langle\nabla\theta,-J(0,D\eta)\rangle d\mu_{g}\\
&  =\int_{\Omega}\langle\nabla\theta,(D\eta,0)\rangle d\mu_{g}.
\end{align*}
with all inner products thus far being computed with respect to the ambient
Euclidean metric. Now
\[
\nabla\theta=g^{ij}\theta_{i}\partial_{j}%
\]
where
\begin{align*}
\partial_{1}  &  =(1,0,\dots,0,u_{11},u_{21},\dots,u_{n1}),\\
&  ...,\\
\partial_{n}  &  =(0,0,\dots,1,u_{1n},u_{2n},\dots,u_{nn}),
\end{align*}
so we have%
\begin{align*}
\delta F_{\Omega}(\eta)  &  =\int_{\Omega}\left\langle g^{ij}\theta
_{i}\partial_{j},(D\eta,0)\right\rangle d\mu_{g}\\
&  =\int_{\Omega}g^{ij}\theta_{i}\eta_{j}d\mu_{g}\\
&  =\int_{\Omega}\left\langle \nabla\theta,\nabla\eta\right\rangle _{g}%
d\mu_{g}.
\end{align*}
Thus we have
\[
\delta F_{\Omega}(\eta)=0\ \ \ \mbox{ for all $\eta \in C_{0}^{\infty
}(\Omega)$}
\]
if and only if
\[
\int_{\Omega}\left\langle \nabla\theta,\nabla\eta\right\rangle d\mu
_{g}=0\ \ \text{ for all }\eta\in C_{c}^{\infty}(\Omega).
\]
This equation has the weak form
\[
\int_{\Omega}\eta\Delta_{g}\theta d\mu_{g}=0\text{ for all }\eta\in
C_{c}^{\infty}(\Omega)
\]
that is
\begin{equation}
\Delta_{g}\theta=0.
\end{equation}
It follows that for $u\in C^{3}(\Omega)$, the volume (\ref{V2}) is stationary
under Hamiltonian variations precisely when (\ref{Hamstat}) is satisfied.
\ Because (\ref{V1}) and (\ref{V2}) are the same functional, if follows that
for $u\in C^{3}(\Omega),$ (\ref{vHamstat}) and (\ref{Hamstat}) are equivalent. \ 
\end{proof}

Observe that, for the gradient graph $\Gamma_{u}=\left\{  (x,Du(x)):x\in
\Omega\right\}  $, the vertical variations constructed by (\ref{verteta}) are
in 1-1 correspondence with $C_{c}^{\infty}(\Omega)$. Note that one can also
construct a variational field, $V=J\nabla\eta$ for each $\eta\in C_{c}%
^{\infty}(\Gamma_{u}).$ \ This is the traditional way of producing Hamiltonian
variations along any Lagrangian submanifold, graphical or not. If the
potential $u$ is smooth, then $C_{c}^{\infty}(\Gamma_{u})=$ $C_{c}^{\infty
}(\Omega)$ and the sets of variations are in 1-1 correspondence. One can then
compute geometrically
\begin{align}
\left.  \frac{d}{dt}F_{\Omega}(\gamma(t)) \right|  _{t=0}  &  =\int_{\Omega
}\left\langle -\vec{H},V\right\rangle d\mu_{g}\label{vg}\\
&  =\int_{\Omega} \left\langle J\nabla\theta,J\nabla\eta\right\rangle d\mu
_{g}\nonumber\\
&  =\int_{\Omega} \left\langle \nabla\theta,\nabla\eta\right\rangle d\mu_{g}
.\nonumber
\end{align}
In particular, the first variational formula is the same.

Note that, in general, when $u$ is not smooth, we have $C_{c}^{\infty}%
(\Gamma_{u})\neq C_{c}^{\infty}(\Omega).$ \ For example if the submanifold
$\Gamma_{u}$ is smooth but the gradient graph has vertical tangents, one would
expect nearby Lagrangian manifolds that are not graphical: These clearly
cannot be reached through a path of vertical variations. In this case, we have
strict containment%
\[
C_{c}^{\infty}(\Omega)\subsetneq C_{c}^{\infty}(\Gamma_{\Omega}).
\]
Thus a Hamiltonian stationary manifold whose volume is stationary under the
larger set of variations, satisfies the equation (\ref{vHamstat}) as well.
Thus in this case, (\ref{vHamstat}) is formally weaker than (\ref{Hamstat}).
It is worth asking when these equations are the same: We delve into this in
the next section.

We note, as it will become useful later, that if $D^{2}u$ is bounded by a
fixed constant almost everywhere, then from (\ref{metricgij}) we see that the
operator
\[
\Delta_{g}=\frac{1}{\sqrt{g}}\partial_{i}(\sqrt{g}g^{ij}\partial_{j})
\]
is uniformly elliptic.

\section{Proof of Theorem \ref{T1}}

First we will consider a general fourth order Euler-Lagrange type equation of
the form
\begin{equation}
\int a^{ijkl}(D^{2}u)u_{ik}\eta_{jl}dx=0 \label{var}%
\end{equation}
for all $\eta\in C_{c}^{\infty}$, where each $a^{ijkl}$ is a smooth function
defined on Hessian space. A function $u\in W^{2,\infty}(\Omega) $ is called a
variational solution to \eqref{var} on $\Omega$, if (\ref{var}) is satisfied
for all $\eta\in C_{c}^{\infty}(\Omega).$ (The choice of the space
$W^{2,\infty}( \Omega) $ may not be the most general, however, it suffices for
our purposes since we will only be considering the case when $u\in C^{1,1}.$)

The proof of the following lemma is based on the calculation in \cite[section
6.3]{Evans}. Essentially, if we have a fourth order nonlinear elliptic
equation\ of type (\ref{var}) such that the nonlinearity $a^{ijkl}(D^{2}u)$
has either a mild or `monotone' dependence on $D^{2}u,$ we can prove increased
regularity for solutions of the equation.

\begin{lem}
\label{IR2} \bigskip Suppose that $u\in W^{2,\infty}\left(  \Omega\right)  $
is a weak solution to (\ref{var}) on $\Omega$ for $n\geq2.$ Suppose there is a
convex neighborhood in Hessian space $U\subset S^{n\times n\text{ }}$such that
for all $M,M^{\ast},M^{\prime}\in U$
\begin{equation}
\frac{\partial a^{ijkl}}{\partial u_{pq}}(M^{\ast})M_{ik}^{\prime}W_{pq}%
W_{jl}+a^{ijkl}(M)W_{ik}W_{jl}\geq\beta\sum_{r,s}W_{rs}^{2} \label{condition4}%
\end{equation}
for all symmetric matrices $W$, where $\beta$ is a positive constant. If
$D^{2} u \left(  \Omega\right)  \subset U,$ wherever $D^{2}\bigskip u$ is
defined, then $u\in W_{loc}^{3,2}(\Omega)$.
\end{lem}

\begin{proof}
By approximation, the equation \eqref{var} must hold for compactly supported
test functions in $W_{0}^{2,\infty}\left(  \Omega\right)  $; in particular, it
must hold for the double difference quotient
\[
\eta=-\left[  \zeta^{4}u^{\left(  h_{m}\right)  }\right]  ^{\left(
-h_{m}\right)  }%
\]
where $\zeta$ $\in C_{c}^{\infty}\left(  \Omega\right)  $ is a cutoff function
that is $1$ on some interior set, and the upper $\left(  h_{m}\right)  $
refers to the difference quotient%
\[
f^{\left(  h_{m}\right)  }(x):=\frac{f\left(  x+he_{m}\right)  -f(x)}{h}%
\]
and we have chosen $h$ small enough (depending on $\zeta$) so that $\eta$ is
well defined and compactly supported. We have%
\begin{equation}
\int_{\Omega}a^{ijkl}(D^{2}u)u_{ik}\left(  -\left[  \zeta^{4}u^{\left(
h_{m}\right)  }\right]  ^{\left(  -h_{m}\right)  }\right)  _{jl}dx=0.
\label{EL31}%
\end{equation}
For $h$ small, we can \textquotedblleft integrate by parts" with respect to
the difference quotient, i.e.
\[
\int_{\Omega}\left[  a^{ijkl}(D^{2}u)u_{ik}\right]  ^{\left(  h_{m}\right)
}\left(  \zeta^{4}u^{\left(  h_{m}\right)  }\right)  _{jl}dx=0.
\]
Now the \textquotedblleft product rule" for difference quotients gives%
\begin{align*}
\left[  a^{ijkl}(D^{2}u)u_{ik}\right]  ^{\left(  h_{m}\right)  }(x)  &
=u_{ik}(x+he_{m})\frac{a^{ijkl}(D^{2}u(x+he_{m}))-a^{ijkl}(D^{2}u(x))}{h}\\
&  \ \ \ +a^{ijkl}(D^{2}u(x))\frac{u_{ik}(x+he_{m})-u_{ik}(x)}{h}\\
&  =u_{ik}(x+he_{m})\int_{0}^{1}\frac{\partial a^{ijkl}}{\partial u_{pq}%
}\left(  (1-t)D^{2}u(x)+tD^{2}u(x+he_{m})\right)  \frac{u_{pq}(x+he_{m}%
)-u_{pq}(x)}{h}dt\\
&  \ \ \ +a^{ijkl}(D^{2}u(x))\frac{u_{ik}(x+he_{m})-u_{ik}(x)}{h}\\
&  =A^{ijkl,pq}(x)u_{ik}(x+he_{m})v_{pq}(x)+a^{ijkl}(D^{2}u(x))v_{ik}(x)
\end{align*}
where
\[
v=u^{\left(  h_{m}\right)  }%
\]
and
\begin{align*}
A^{ijkl,pq}(x)  &  =\int_{0}^{1}\frac{\partial a^{ijkl}}{\partial u_{pq}%
}\left(  (1-t)D^{2}u(x)+tD^{2}u(x+he_{m})\right)  dt\\
&  =\frac{\partial a^{ijkl}}{\partial u_{pq}}\left(  M^{\ast}(x)\right)
\end{align*}
where
\[
M^{\ast}(x):=(1-t^{\ast})D^{2}u(x)+t^{\ast}D^{2}u(x+he_{m})
\]
for some $t^{\ast}$ by the mean value theorem. (Note that for a fixed $h$,
$D^{2}u$ exists at both $x$ and $x+he_{m},$ almost everywhere, so all of the
above quantities are defined almost everywhere.) So equation (\ref{EL31})
becomes%
\[
\int_{\Omega}\left(  \frac{\partial a^{ijkl}}{\partial u_{pq}}\left(  M^{\ast
}(x)\right)  u_{ik}(x+he_{m})v_{pq}(x)+a^{ijkl}(D^{2}u(x))v_{ik}(x)\right)
\left(  \zeta^{4}v(x)\right)  _{jl}dx=0.
\]
Now differentiating the second factor,%
\begin{equation}
\int_{\Omega}\left(
\begin{array}
[c]{c}%
\left(  \frac{\partial a^{ijkl}}{\partial u_{pq}}\left(  M^{\ast}(x)\right)
u_{ik}(x+he_{m})v_{pq}(x)+a^{ijkl}(D^{2}u(x))v_{ik}(x)\right) \\
\times\left[  \zeta^{4}v_{jl}+4\zeta^{3}\zeta_{j}v_{l}+4\zeta^{3}\zeta
_{l}v_{j}+4v(\zeta^{3}\zeta_{jl}+3\zeta^{2}\zeta_{j}\zeta_{l})\right]  (x)
\end{array}
\right)  dx=0. \label{bigc1}%
\end{equation}

By the condition (\ref{condition4})\ in the hypothesis we have that%
\[
\int_{\Omega}\left(  \frac{\partial a^{ijkl}}{\partial u_{pq}}\left(  M^{\ast
}(x)\right)  u_{ik}(x+he_{m})v_{pq}(x)+a^{ijkl}(D^{2}u(x))v_{ik}(x)\right)
\zeta^{4}v_{jl}dx\geq\beta\int_{\Omega}\zeta^{4}\sum_{r,s}v_{rs}^{2}dx.
\]

For the remaining terms, note that for the second term in the expansion of
(\ref{bigc1}) we have by Young's inequality
\begin{align*}
&  \left\vert \frac{\partial a^{ijkl}}{\partial u_{pq}}\left(  M^{\ast
}(x)\right)  u_{ik}(x+he_{m})v_{pq}(x)4\zeta^{3}(x)\zeta_{j}(x)v_{l}%
(x)\right\vert \leq\\
&  C(n)\frac{1}{\varepsilon}\left(  \frac{\partial a^{ijkl}}{\partial u_{pq}%
}\left(  M^{\ast}(x)\right)  \right)  ^{2}\left(  u_{ik}(x+he_{m})\right)
^{2}\zeta^{2}(x)\left\vert D\zeta(x)\right\vert ^{2}|Dv(x)|^{2}+\varepsilon
\zeta^{4}(x)v_{pq}^{2}(x).
\end{align*}
A similar expression can be made for each of the terms. Noting that $D^{2}u$
is bounded and $v$ is the different quotient of $u$, we obtain%
\begin{align*}
&  \int_{\Omega}\left(
\begin{array}
[c]{c}%
\frac{\partial a^{ijkl}}{\partial u_{pq}}\left(  M^{\ast}(x)\right)
u_{ik}(x+he_{m})v_{pq}(x)+a^{ijkl}(D^{2}u(x))v_{ik}(x)\\
\times\left[  4\zeta^{3}\zeta_{j}v_{l}+4\zeta^{3}\zeta_{l}v_{j}+4v(\zeta
^{3}\zeta_{jl}+3\zeta^{2}\zeta_{j}\zeta_{l})\right]  (x)
\end{array}
\right)  dx\\
&  \leq C(|D^{{}}u|,|D^{2}u|,\left\vert D\zeta\right\vert ,\left\vert
D^{2}\zeta\right\vert ^{2},\left\vert Da^{ijkl}\right\vert )\frac
{1}{\varepsilon}\int_{\Omega}|Dv|^{2}dx+\varepsilon\int_{\Omega}\sum
_{r,s}\zeta^{4}v_{rs}^{2}dx
\end{align*}
where $\left\vert Da^{ijkl}\right\vert $ is a norm on the total derivative of
the functions $a^{ijkl}$ on the space of symmetric matrices.

We conclude that by choosing $\varepsilon$ appropriately, we have
\begin{align*}
\frac{\beta}{2}\int_{\Omega}\zeta^{4}\sum_{r,s}v_{rs}^{2}dx  &  \leq C(|Du|,
|D^{2}u|,\left\vert D\zeta\right\vert ,\left\vert D^{2}\zeta\right\vert
^{2},\left\vert Da^{ijkl}\right\vert )\frac{1}{\varepsilon}\int_{\Omega
}|Dv|^{2}dx\\
&  \leq C\left\Vert v\right\Vert _{W^{1,2}(\Omega)}\\
&  \leq C\left\Vert u\right\Vert _{W^{2,2}(\Omega)}.
\end{align*}
Thus
\[
\left\Vert v\right\Vert _{W^{2,2}(\left\{  x|\zeta(x)=1\right\}  )}\leq C.
\]
Now this estimate is uniform in $h$ and direction $e_{m}$ so we conclude that
the derivatives are in $W^{2,2\text{ }}(\Omega)$ and thus $u\in W^{3,2}%
(\left\{  x|\zeta(x)=1\right\}  )$.
\end{proof}

\begin{prop}
\label{IR32} There is a bound $c(n)$ such that if
\[
\left\Vert u\right\Vert _{C^{1,1}\left(  \Omega\right)  }\leq c(n)
\]
for a weak solution $u$ to the Hamiltonian stationary equation (\ref{vHamstat}%
), then $u\in W_{loc}^{3,2}(\Omega)$ . \ 
\end{prop}

\begin{proof}
First recall (cf. \cite[section 5.8.2]{Evans}) that the Hessian $D^{2}u$ is
defined almost everywhere and bounded where it is defined in terms of the
$C^{1,1}$ norm. Considering (\ref{vHamstat}) in the notation of (\ref{var}) we
have
\[
a^{ijkl}=\sqrt{g}g^{ij}\delta^{kl}.
\]
Our goal is to show that the condition (\ref{condition4}) is satisfied on the
set
\[
U=\left\{  M\in S^{n\times n}:\Vert M\Vert_{\infty}\leq c(n)\right\}  .
\]
For simplicity, we will write $|M|$ for $\Vert M\Vert_{\infty}$, especially
when Hessian is involved.

Computing, we see
\begin{align}
\frac{\partial a^{ijkl}}{\partial u_{mp}}  &  =\frac{1}{2}\sqrt{g}g^{ab}%
\frac{\partial}{\partial u_{mp}}g_{ab}g^{ij}\delta^{kl}-\sqrt{g}g^{ia}%
g^{bj}\frac{\partial}{\partial u_{mp}}g_{ab}\delta^{kl}\label{comp1}\\
&  =\left(  \frac{1}{2}g^{ab}g^{ij}\delta^{kl}-g^{ia}g^{bj}\delta^{kl}\right)
\sqrt{g}\frac{\partial}{\partial u_{mp}}g_{ab}\nonumber\\
&  =\left(  \frac{1}{2}g^{ab}g^{ij}\delta^{kl}-g^{ia}g^{bj}\delta^{kl}\right)
\sqrt{g}\frac{\partial}{\partial u_{mp}}\left(  \delta_{ab}+u_{ac}\delta
^{cd}u_{db}\right) \nonumber\\
&  =\left(  \frac{1}{2}g^{ab}g^{ij}\delta^{kl}-g^{ia}g^{bj}\delta^{kl}\right)
\sqrt{g}\left(  \delta_{mp,ac}\delta^{cd}u_{db}+u_{ac}\delta^{cd}%
\delta_{mp,db}\right)  .\nonumber
\end{align}
In particular,%
\begin{equation}
\left\vert \frac{\partial a^{ijkl}}{\partial u_{pq}}(D^{2}u)\right\vert \leq
C(n)\left\vert D^{2}u\right\vert \left(  1+\left\vert D^{2}u\right\vert
^{2}\right)  ^{n/2}. \label{nlpert}%
\end{equation}
Next, note that if we let
\[
G_{ij}=\sqrt{g}g^{ij},
\]
we can write%
\[
\sqrt{g}g^{ij}\delta^{kl}W_{ik}W_{jl}=\mbox{Trace}(G^{T}WI_{n}W^{T}).
\]
But $G$ can be diagonalized by an orthogonal matrix $O:$
\[
G^{T}=O^{T}DO
\]
where
\[
D=\sqrt{g}\left(
\begin{array}
[c]{ccc}%
\frac{1}{1+\lambda_{1}^{2}} & 0 & 0\\
0 & ... & 0\\
0 & 0 & \frac{1}{1+\lambda_{n}^{2}}%
\end{array}
\right)  .
\]
Then
\begin{align*}
\sqrt{g}g^{ij}\delta^{kl}W_{ik}W_{jl}  &  ={\mbox{Trace}}(O^{T}DOWW^{T})\\
&  ={\mbox{Trace}}(OO^{T}DOWW^{T}O^{T})\\
&  ={\mbox{Trace}}(D\left(  OW\right)  (OW)^{T})\\
&  \geq\min_{i}D_{ii}\cdot\mbox{Trace}\left(  (OW\right)  (OW)^{T})\\
&  =\min_{i}D_{ii}\left\Vert OW\right\Vert _{HS}^{2}\\
&  =\min_{i}D_{ii}\left\Vert W\right\Vert _{HS}^{2},
\end{align*}
where we are using the Hilbert-Schmidt norm on matrices. Thus
\begin{equation}
\sqrt{g}g^{ij}\delta^{kl}W_{ik}W_{jl}\geq\frac{1}{1+c(n)^{2}}\left\Vert
W\right\Vert _{HS}^{2}. \label{elli5}%
\end{equation}
Combining (\ref{nlpert}) and (\ref{elli5}) and plugging this into
(\ref{condition4}) we see for $M^{\ast},M^{\prime},$ and $M$ in $U$ we have
\begin{align*}
&  \frac{\partial a^{ijkl}}{\partial u_{pq}}(M^{\ast})M_{ik}^{\prime}%
W_{pq}W_{jl}+a^{ijkl}(M)W_{ik}W_{jl}\\
&  \geq\frac{1}{1+c(n)^{2}}\left\Vert W\right\Vert _{HS}^{2}-C(n)\left\vert
c(n)\right\vert ^{2}\left(  1+c(n)^{2}\right)  ^{n/2}\left\Vert W\right\Vert
_{\infty}^{2}\\
&  \geq\beta\left\Vert W\right\Vert _{HS}^{2}%
\end{align*}
for some $\beta>0$, using the equivalence of norms, when $c(n)$ is chosen
small. \ The conclusion follows from Lemma \ref{IR2}.
\end{proof}

To extend solutions across a small set in Theorem \ref{T1}. we will need the
following theorem of Serrin (Theorem 2 in \cite{Serrin}).

\begin{theorem}
\label{Serrin}(Serrin) Suppose $n\geq2$ and that $f$ is a\ bounded continuous
weak solution to a uniformly elliptic second order divergence equation with
bounded measurable coefficients on ${\Omega}-Q,$ for an open domain ${\Omega}$
and $Q$ a compact subset. If $Q$ has capacity zero, then $f$ may be extended
to a weak solution across the domain ${\Omega}.$
\end{theorem}

We now proceed to prove Theorem \ref{T1}.

\begin{proof}
First, let us consider the case when $Q$ is the empty set. Because $u\in
W_{loc}^{3,2}\left(  \Omega\right)  \cap C^{1,1}(\Omega)$ we may use a
standard mollification construction, letting
\[
u^{\varepsilon}=\rho_{\varepsilon}\ast u
\]
for an appropriate function $\rho_{\varepsilon}$ as in \cite[Appendix
C.4]{Evans}. In particular (see \cite[Appendix C, Theorem 6]{Evans})
\[
\lim_{\varepsilon\rightarrow0}\left\Vert u^{\varepsilon}-u\right\Vert
_{W_{loc}^{3,2}\left(  \Omega\right)  }=0
\]
and each $u^{\varepsilon}$ is smooth.

Now we define functionals on $C_{c}^{\infty}\left(  \Omega\right)  $ by
\begin{align*}
F^{\varepsilon}(\eta)  &  =\int_{\Omega}\left[  \sqrt{g}g^{ij}\delta
^{kl}u_{ik}\right]  ^{\varepsilon}\eta_{jl}dx\\
F(\eta)  &  =\int_{\Omega}\sqrt{g}g^{ij}\delta^{kl}u_{ik}\eta_{jl}dx
\end{align*}
with the notation $\left[  \sqrt{g}g^{ij}\delta^{kl}u_{ik}\right]
^{\varepsilon}$ meaning \textquotedblleft constructed from $u^{\varepsilon}$
using (\ref{metricgij}) ,$"$ (in particular, this does \emph{not} mean the
mollification of the expression). \ 

First we check that for each $\eta,$%
\[
F(\eta)=\lim_{\varepsilon\rightarrow0}F^{\varepsilon}(\eta).
\]
We have
\begin{align*}
F^{\varepsilon}(\eta)-F(\eta)  &  =\int_{\Omega}\left(  \left[  \sqrt{g}%
g^{ij}u_{ik}\right]  ^{\varepsilon}-\sqrt{g}g^{ij}u_{ik}\right)  \delta
^{kl}\eta_{jl}dx\\
&  =\int_{\Omega}\left(  \left[  \sqrt{g}g^{ij}u_{ik}\right]  ^{\varepsilon
}-\left[  \sqrt{g}g^{ij}\right]  ^{\varepsilon}u_{ik}+\left[  \sqrt{g}%
g^{ij}\right]  ^{\varepsilon}u_{ik}-\sqrt{g}g^{ij}u_{ik}\right)  \delta
^{kl}\eta_{jl}dx\\
&  =\int_{\Omega}\left(  \left[  \sqrt{g}g^{ij}\right]  ^{\varepsilon}\left(
u_{ik}^{\varepsilon}-u_{ik}\right)  +\left(  \left[  \sqrt{g}g^{ij}\right]
^{\varepsilon}-\sqrt{g}g^{ij}\right)  g^{ij}u_{ik}\right)  \delta^{kl}%
\eta_{jl}dx
\end{align*}
Now because $u\in C^{1,1}$ and $\eta_{jl}$ is bounded, we simply have to check
that
\begin{align*}
\left(  u_{ik}^{\varepsilon}-u_{ik}\right)   &  \rightarrow0\ \text{\ in
}L_{loc}^{1}\\
\left(  \left[  \sqrt{g}g^{ij}\right]  ^{\varepsilon}-\sqrt{g}g^{ij}\right)
&  \rightarrow0\text{\ in }L_{loc}^{1}.
\end{align*}
The first assertion is clear as $u\in W_{loc}^{3,2}\left(  \Omega\right)  .$

Next,
\[
\left\vert \left[  \sqrt{g}g^{ij}\right]  ^{\varepsilon}-\sqrt{g}
g^{ij}\right\vert \leq\sup_{i,j}\left\vert \frac{\partial\left(  \sqrt
{g}g^{ij}\right)  }{\partial u_{ab}}\right\vert \left(  u_{ab}^{\varepsilon
}-u_{ab}\right)  .
\]
Mimicking computations following (\ref{comp1}) we see
\[
\left\vert \frac{\partial\left(  \sqrt{g}g^{ij}\right)  }{\partial u_{ab}
}\right\vert \leq C(n)\left\vert D^{2}u\right\vert \left(  1+\left\vert
D^{2}u\right\vert ^{2}\right)  ^{n/2}\leq C.
\]
Thus
\begin{equation}
\left\vert \left[  \sqrt{g}g^{ij}\right]  ^{\varepsilon}-\sqrt{g}
g^{ij}\right\vert \leq C\left\vert D^{2}u^{\varepsilon}-D^{2}u\right\vert
\label{bound_2}%
\end{equation}
and the second assertion then follows from the first.

We conclude that
\[
F(\eta)=\lim_{\varepsilon\rightarrow0}F^{\varepsilon}(\eta).
\]

Next, we define functionals
\begin{align*}
G^{\varepsilon}(\eta)  &  =\int_{\Omega}\left[  \sqrt{g}g^{ij}\theta
_{i}\right]  ^{\varepsilon}\eta_{j}dx\\
G(\eta)  &  =\int_{\Omega}\sqrt{g}g^{ij}\theta_{i}\eta_{j}dx=\int_{\Omega
}\sqrt{g}g^{ij} g^{ab}u_{abi}\eta_{j}dx
\end{align*}
recalling that
\[
\theta_{i}=\left(  \operatorname{Im}\log\det\left(  I+iD^{2}u\right)  \right)
_{i}=g^{ab}u_{abi}
\]
and noting that since $u\in W_{loc}^{3,2}\left(  \Omega\right)  $, the third
derivatives exist almost everywhere.

Applying the first variational formulae for smooth submanifolds in section 2
to the smooth $\Gamma_{u^{\varepsilon}}$, we see that
\[
\delta F_{\Omega}(\eta)=\int_{\Omega}\left[  \sqrt{g}g^{ij}\delta^{kl}%
u_{ik}\right]  ^{\varepsilon}\eta_{jl}dx=\int_{\Omega}\left[  \sqrt{g}%
g^{ij}\theta_{i}\right]  ^{\varepsilon}\eta_{j}dx
\]
that is
\[
G^{\varepsilon}(\eta)=F^{\varepsilon}(\eta).
\]
So clearly, from our observations on $F^{\varepsilon}(\eta)$ we see that
\[
\lim_{\varepsilon\rightarrow0}G^{\varepsilon}(\eta)=0.
\]
All that remains is to show that
\[
\lim_{\varepsilon\rightarrow0}G^{\varepsilon}(\eta)=G(\eta).
\]
We follow the same procedure as above:
\begin{align*}
G^{\varepsilon}(\eta)-G(\eta)  &  =\int_{\Omega}\left(  \left[  \sqrt{g}%
g^{ij}\theta_{i}\right]  ^{\varepsilon}-\sqrt{g}g^{ij}\theta_{i}\right)
\eta_{j}dx\\
&  =\int_{\Omega}\left(  \left[  \sqrt{g}g^{ij}\theta_{i}\right]
^{\varepsilon}-\left[  \sqrt{g}g^{ij}\right]  ^{\varepsilon}\theta_{i}+\left[
\sqrt{g}g^{ij}\right]  ^{\varepsilon}\theta_{i}-\sqrt{g}g^{ij}\theta
_{i}\right)  \eta_{j}dx\\
&  =\int_{\Omega}\left(  \left[  \sqrt{g}g^{ij}\right]  ^{\varepsilon}\left(
\theta_{i}^{\varepsilon}-\theta_{i}\right)  +\left(  \left[  \sqrt{g}%
g^{ij}\right]  ^{\varepsilon}-\sqrt{g}g^{ij}\right)  \theta_{i}\right)
\eta_{j}dx
\end{align*}
Now we have to be slightly more careful, but proceed as before: Starting with
the last term, we have using (\ref{bound_2})
\begin{align*}
\int_{\Omega}\left(  \left[  \sqrt{g}g^{ij}\right]  ^{\varepsilon}-\sqrt
{g}g^{ij}\right)  \theta_{i}\eta_{j}dx  &  \leq\left\Vert D\theta\right\Vert
_{L^{2}}\left\Vert D\eta\right\Vert _{L^{\infty}}\left\Vert \left[  \sqrt
{g}g^{ij}\right]  ^{\varepsilon}-\sqrt{g}g^{ij}\right\Vert _{L^{2}}\\
&  \leq\left\Vert D\theta\right\Vert _{L^{2}}\left\Vert D\eta\right\Vert
_{L^{\infty}}C\left\Vert D^{2}u^{\varepsilon}-D^{2}u\right\Vert _{L^{2}}\\
&  \rightarrow0
\end{align*}
as
\[
\left\Vert D\theta\right\Vert _{L^{2}\left(  K\right)  }\leq C\left\Vert
u\right\Vert _{W^{3,2}(K)}%
\]
for any $K$ compact inside $\Omega.$\ 

Finally
\begin{align*}
&  \int_{\Omega}\left[  \sqrt{g}g^{ij}\right]  ^{\varepsilon}\left(
\theta_{i}^{\varepsilon}-\theta_{i}\right)  \eta_{j}dx\\
&  =\int_{\Omega}\left[  \sqrt{g}g^{ij}\right]  ^{\varepsilon}\left(  \left(
g^{ab}\right)  ^{\varepsilon}u_{abi}^{\varepsilon}-\left(  g^{ab}\right)
^{\varepsilon}u_{abi}+\left(  g^{ab}\right)  ^{\varepsilon}u_{abi}%
-g^{ab}u_{abi}\right)  \eta_{j}dx\\
&  \leq\sup\left[  \sqrt{g}g^{ij}\right]  \left\Vert D\eta\right\Vert
_{L^{\infty}}\left\{  \left\Vert \left(  g^{ab}\right)  ^{\varepsilon
}\right\Vert _{L^{2}}\left\Vert u_{abi}^{\varepsilon}-u_{abi}\right\Vert
_{L^{2}}+\left\Vert \left(  g^{ab}\right)  ^{\varepsilon}-g^{ab}\right\Vert
_{L^{2}}\left\Vert D^{3}u\right\Vert _{L^{2}}\right\}  .
\end{align*}
Because $u^{\varepsilon}\rightarrow u$ in $W_{loc}^{3,2},$ these terms go to zero.

We conclude that
\[
G(\eta)=\int_{\Omega}\sqrt{g}g^{ij}\theta_{i}\eta_{j}dx=0
\]
for all test functions $\eta.$ \ It follows that $\theta$ is a weak solution
of the uniformly elliptic equation (\ref{Hamstat}).


When $Q$ is a compact subset in $\Omega$, because ${\Omega}\backslash Q$ is
itself an open domain, the result established above asserts that $u\in
W_{loc}^{3,2}\left(  {\Omega}\backslash Q\right)  $ and $u$ is a weak solution
to (\ref{Hamstat}) on ${\Omega}\backslash Q.$ This means that
(\ref{weakintegral})\ holds for all $\eta$ supported in ${\Omega}$ away from
$Q$. So $\theta$ is now in the setting of Serrin's Theorem: We can extend
$\theta$ to a weak solution across the entire domain, so $u$ is a weak
solution to (\ref{Hamstat}) on ${\Omega}$. Next, we apply Theorem \ref{T2}
(whose proof is independent of Theorem \ref{T1}), where the condition
\eqref{delta3} applies. We conclude that $u$ is smooth on ${\Omega}$. Thus,
the first variation formulae yield equivalence of (\ref{Hamstat}) and
(\ref{vHamstat}), so $u$ must be a solution of (\ref{vHamstat}) on ${\Omega.}$
\end{proof}

\section{Lewy-Yuan rotations}

In this section we discuss and motivate the Lewy-Yuan rotation. \ We risk
giving extra descriptions here in order to give a clear motivation as to what
the rotation is useful for. \ We also rigorously justify low regularity
versions of the Lewy-Yuan rotation.

In the special Lagrangian setting, Yuan \cite{YYI} used the following unitary
change of coordinates
\begin{align}
U  &  :\mathbb{C}^{n}\rightarrow\mathbb{C}^{n}\label{unitary1}\\
U(x+\sqrt{-1}y)  &  =e^{-\sqrt{-1}\pi/4}\left(  x+\sqrt{-1}y\right)
.\nonumber
\end{align}
In this case, a surface $\Gamma$ that was the gradient graph of a convex
function $u$ over the original $\mathbb{R}^{n}$-plane, is now represented as a
gradient graph of a new function $\bar{u}$ over the new $\mathbb{R}^{n}%
$-plane, but this time with
\[
-I_{n}\leq D^{2}\bar{u}\leq I_{n}.
\]
We call this a downward rotation by angle $\pi/4:$ \ The word `downward'
refers to the fact that the argument of the complex number $e^{-\sqrt{-1}%
\pi/4}$ (\ref{unitary1}) is negative. Any surface $\Gamma$ that is the
gradient graph of a semi-convex function $u$ can be rotated downward
(\cite{YY06}). \ If for $\beta\in(0,\pi/2)$ we have
\[
D^{2}u\geq-\tan\beta\,I_{n}%
\]
then we can rotate the graph downward by any positive angle $\alpha
<\pi/2-\beta.$ \ More precisely, given
\[
\Gamma=\left\{  \left(  x,Du(x)\right)  ,x\in\Omega\right\}  \subset
\mathbb{R}^{n}+\sqrt{-1}\mathbb{R}^{n}%
\]
over $\Omega$, let
\begin{equation}
\bar{\Gamma}=U_{\alpha}\Gamma\label{rotation}%
\end{equation}
where%
\begin{equation}
U_{\alpha}=\left(
\begin{array}
[c]{ccc}%
e^{-\sqrt{-1}\alpha} &  & \\
& ... & \\
&  & e^{-\sqrt{-1}\alpha}%
\end{array}
\right)  . \label{unitary}%
\end{equation}
Clearly, $\bar{\Gamma}$ is isometric to $\Gamma$ via the unitary rotation. In
coordinates, this is equivalent to the following map.
\begin{align}
\bar{x}  &  =\cos(\alpha)x+\sin(\alpha)Du(x)\label{xbarybar}\\
\bar{y}  &  =-\sin(\alpha)x+\cos(\alpha)Du(x).\nonumber
\end{align}
Here $\bar{x}$ and $\bar{y}$ are simply the projections onto $\mathbb{R}^{n}$
and $\sqrt{-1}\mathbb{R}^{n}$ of $\bar{\Gamma}$, respectively.

Considering the functions $\bar{x}(x),\bar{y}(x)$ we may compute the
differential form
\begin{align*}
\sum_{i}\bar{y}^{i}d\bar{x}^{i}  &  =\sum_{i}\left(  -\sin(\alpha)x^{i}%
+\cos(\alpha)u_{i}(x)\right)  \left(  \cos(\alpha)dx^{i}+\sin(\alpha
)u_{ij}(x)dx^{j}\right) \\
&  =\sum_{i}\left(
\begin{array}
[c]{c}%
-\sin(\alpha)\cos(\alpha)x^{i}dx^{i}+\cos^{2}(\alpha)u_{i}(x)dx^{i}\\
-\sin^{2}(\alpha)x^{i}u_{ij}(x)dx^{j}+\cos(\alpha)\sin(\alpha)u_{i}%
(x)u_{ij}(x)dx^{j}%
\end{array}
\right) \\
&  =-\sin(\alpha)\cos(\alpha)D\frac{\left\vert x\right\vert ^{2}}{2}+\cos
^{2}(\alpha)Du(x)\\
&  \ \ \ \ -\sin^{2}(\alpha)\left(  D(x\cdot Du(x)\right)  -Du(x))+\cos
(\alpha)\sin(\alpha)D\frac{\left\vert Du(x)\right\vert ^{2}}{2}\\
&  =Du+\sin(\alpha)\cos(\alpha)D\frac{\left\vert Du(x)\right\vert
^{2}-\left\vert x\right\vert ^{2}}{2}-\sin^{2}(\alpha)\left(  D(x\cdot
Du\right)  )\\
&  =D\left(  u(x)+\sin(\alpha)\cos(\alpha)\frac{\left\vert Du(x)\right\vert
^{2}-\left\vert x\right\vert ^{2}}{2}-\sin^{2}(\alpha)\left(  (x\cdot
Du(x)\right)  )\right)  .
\end{align*}
We see that the 1-form $\sum_{i}\bar{y}^{i}d\bar{x}^{i}$ is exact (regardless
of cohomological conditions) as we can exhibit $\bar{u}\left(  \bar{x}\right)
=\bar{u}\left(  \bar{x}(x)\right)  $ solving $D_{\bar{x}}\bar{u}%
=\allowbreak\bar{y}d\bar{x}^{i}.$ It follows that
\[
(\bar{x},\bar{y})=(\bar{x},D_{\bar{x}}\bar{u}(\bar{x}))
\]
for some function $\bar{u}\left(  \bar{x}\right)  $. The potential $\bar{u}$
is given explicitly, however, the explicit formula is only given in terms of
the $x$ coordinates. Fortunately, $\bar{x}(x)$ is a change of coordinates
(this follows from the semi-convexity, see Proposition \ref{weakrotation}
below) and is invertible.

To summarize, we have exhibited $\bar{\Gamma}$ both as the gradient graph of a
function $\bar{u}$ and as an isometric image of $\Gamma$. The result will be a
new graph with a potential whose Hessian satisfies (see \cite[(1.5) and
(1.6)]{Warren2016})%
\[
-\tan(\beta+\alpha)I_{n}\leq D^{2}\bar{u}\leq\tan(\pi/2-\alpha)I_{n}.
\]
The takeaway is that any semi-convexity guarantees that the graph has a
representation of bounded geometry. Also note that there is nothing sacred
about downward rotations: A function with a Hessian upper bound may always be
rotated upwards to obtain a representation with a Hessian lower bound as well.

Geometrically, if we are not given a potential function, we can always choose
a tangent plane at a point. This plane is Lagrangian, and locally, by the
Poincar\'{e} Lemma, the Lagrangian surface will be a gradient graph over this
tangent plane. In general, one could choose from a large set of unitary
rotations to obtain representations, however, we focus only on the
\textquotedblleft uniform diagonal" rotations of the form (\ref{unitary}) that
rotate each $x$-$y$ plane in the same way.

\subsection{When $\Gamma$ is not smooth}

In the above computation, we referenced the second derivatives of $u,$ despite
the fact that the rotation itself is actually a map on first derivatives.
\ Our goal in this section is to rigorously show that the Lewy-Yuan rotation
can be performed in some low regularity settings where the second derivatives
need not exist everywhere, as long as some semi-convexity is satisfied.

For a constant $K\in\mathbb{R}$, we say that $u$ is $K$-convex on $\Omega$ if
\[
u(x)-K\frac{\left\vert x\right\vert ^{2}}{2}\text{ is convex. }%
\]
For $u$ $\in C^{1}$ this is equivalent to the condition that, for all
$x_{0},x_{1}\in\Omega\,\ $
\begin{equation}
\langle Du(x_{1})-Du(x_{0}),x_{1}-x_{0}\rangle\geq K\left\vert x_{1}%
-x_{0}\right\vert ^{2}. \label{Kconvex}%
\end{equation}

\begin{prop}
\label{weakrotation}Suppose that $\Gamma=\left(  x,Du(x)\right)  $ is a
Lagrangian graph in $\Omega+\sqrt{-1}\mathbb{R}^{n}\subset\mathbb{C}^{n}$ with
$Du$ continuous. \ Suppose that
\begin{equation}
u+\left(  \cot(\sigma)-\varepsilon\right)  \frac{\left\vert x\right\vert ^{2}%
}{2}\text{ \ is convex} \label{semiconvexity}%
\end{equation}
for some $\varepsilon>0,\sigma>0$. Consider the function
\[
\bar{u}(x)=u(x)+\sin\left(  \sigma\right)  \cos\left(  \sigma\right)
\frac{\left\vert Du(x)\right\vert ^{2}-\left\vert x\right\vert ^{2}}{2}%
-\sin^{2}\left(  \sigma\right)  Du(x)\cdot x
\]
and the function $\bar{x}:\Omega\rightarrow\bar{\Omega}\subset\mathbb{R}^{n}$
given by
\begin{equation}
\bar{x}(x)=\cos\left(  \sigma\right)  x+\sin\left(  \sigma\right)  Du(x).
\label{cc}%
\end{equation}
Then

\begin{enumerate}
\item The coordinate change (\ref{cc}) is invertible with Lipschitz continuous inverse,

\item The derivative of $\bar{u}$ in $\bar{x}$ coordinates $\frac{D\bar{u}%
}{d\bar{x}}$ exists everywhere, and

\item The gradient graph $\bar{\Gamma}=\left(  \bar{x},D\bar{u}(\bar
{x}\right)  )\subset\bar{\Omega}+\sqrt{-1}\mathbb{R}^{n}\subset\mathbb{C}^{n}$
is the isometric image of $\ \Gamma$ under the rotation through $\sigma$ as in
(\ref{rotation}).
\end{enumerate}
\end{prop}

\begin{proof}
Note that the convexity condition can be written as, for any two points
$x_{0},x_{1}\in\Omega$,
\[
\left\langle Du(x_{1})-Du(x_{0})+\left(  \cot(\sigma)-\varepsilon\right)
\left(  x_{1}-x_{0}\right)  ,x_{1}-x_{0}\right\rangle \geq0.
\]
This leads to
\begin{equation}
\left\langle \frac{Du(x_{1})-Du(x_{0})}{\left\vert x_{1}-x_{0}\right\vert
},\frac{x_{1}-x_{0}}{\left\vert x_{1}-x_{0}\right\vert }\right\rangle
\geq-\cot(\sigma)+\varepsilon. \label{econvexity}%
\end{equation}
It then follows, for $x_{1}\neq x_{0},$ that
\begin{align}
\left\vert \frac{\bar{x}(x_{1})-\bar{x}(x_{0})}{\left\vert x_{1}%
-x_{0}\right\vert }\right\vert  &  \geq\left\langle \frac{\bar{x}(x_{1}%
)-\bar{x}(x_{0})}{\left\vert x_{1}-x_{0}\right\vert },\frac{x_{1}-x_{0}%
}{\left\vert x_{1}-x_{0}\right\vert }\right\rangle \label{Lips}\\
&  =\left\langle \frac{\cos\left(  \sigma\right)  \left(  x_{1}-x_{0}\right)
+\sin\left(  \sigma\right)  \left(  Du(x_{1})-Du(x_{0})\right)  }{\left\vert
x_{1}-x_{0}\right\vert },\frac{x_{1}-x_{0}}{\left\vert x_{1}-x_{0}\right\vert
}\right\rangle \nonumber\\
&  =\cos\left(  \sigma\right)  +\sin\left(  \sigma\right)  \left\langle
\frac{Du(x_{1})-Du(x_{0})}{\left\vert x_{1}-x_{0}\right\vert },\frac
{x_{1}-x_{0}}{\left\vert x_{1}-x_{0}\right\vert }\right\rangle \nonumber\\
&  \geq\cos\left(  \sigma\right)  -\cot(\sigma)\sin\left(  \sigma\right)
+\sin\left(  \sigma\right)  \varepsilon\nonumber\\
&  =\sin\left(  \sigma\right)  \varepsilon\nonumber
\end{align}
using (\ref{econvexity}). Therefore the continuous map $\bar{x}$ is invertible
and its inverse is Lipschitz continuous with a Lipschitz constant $1/\left(
\sin(\sigma)\varepsilon\right)  $.


Next, for the gradient of $\bar{u}$ in terms of $\bar{x}$, we will compute a
difference quotient%
\[
\bar{u}_{\bar{j}}(\bar{x}_{0})=\lim_{h\rightarrow0}\frac{\bar{u}(\bar{x}%
_{0}+h\bar{e}_{j})-\bar{u}(\bar{x}_{0})}{h}.
\]
Since $\bar{x}$ is invertible, for $\bar{x}_{0}\in\bar{\Omega}$ we may solve,
for small fixed $h$%
\begin{align*}
\bar{x}(x_{0})  &  =\bar{x}_{0}\\
\bar{x}(x_{h})  &  =\bar{x}_{0}+h\bar{e}_{j}%
\end{align*}
that is
\begin{align*}
\cos\left(  \sigma\right)  x_{0}+\sin\left(  \sigma\right)  Du(x_{0})  &
=\bar{x}_{0}\\
\cos\left(  \sigma\right)  x_{h}+\sin\left(  \sigma\right)  Du(x_{h})  &
=\bar{x}_{h}=\bar{x}_{0}+h\bar{e}_{j}.
\end{align*}
Let
\[
\vec{v}=x_{h}-x_{0}.
\]
Then $\vec{v}$ will satisfy%
\begin{equation}
\cos\left(  \sigma\right)  \vec{v}+\sin\left(  \sigma\right)  \left[
Du(x_{h})-Du(x_{0})\right]  =h\bar{e}_{j}. \label{differenceh}%
\end{equation}
Let
\[
\vec{v}=h\vec{V}.
\]
Observe that
\[
\left\vert \vec{V}\right\vert =\frac{\left\vert \vec{v}\right\vert }{h}%
=\frac{\left\vert x_{h}-x_{0}\right\vert }{\left\vert \bar{x}(x_{h})-\bar
{x}(x_{0})\right\vert }\leq\frac{1}{\varepsilon\sin\sigma}%
\]
by (\ref{Lips}). In particular, $\vec{V}$ is a bounded vector. (While the
vector $\vec{V}$ depends on $h,$ we suppress this dependence.) \ The function
$\bar{u}$ is given in term of $x$ coordinates, so in order to evaluate it, we
have to use the change of coordinates, that is%
\[
\bar{u}\left(  \bar{x}_{0}\right)  =\bar{u}(\bar{x}^{-1}(\bar{x}_{0}))=\bar
{u}(x_{0}).
\]
So we may compute the difference quotient of $\bar{u}$ in terms of $x$
\begin{align*}
\frac{\bar{u}\left(  \bar{x}_{h}\right)  -\bar{u}\left(  \bar{x}_{0}\right)
}{h}  &  =\frac{\bar{u}(\bar{x}^{-1}(\bar{x}_{h}))-\bar{u}(\bar{x}^{-1}%
(\bar{x}_{0}))}{h}\\
&  =\frac{u(x_{h})-u(x_{0})}{h}+\sin\left(  \sigma\right)  \cos\left(
\sigma\right)  \frac{\left\vert Du(x_{h})\right\vert ^{2}-\left\vert
Du(x_{0})\right\vert ^{2}-\left\vert x_{h}\right\vert ^{2}+\left\vert
x_{0}\right\vert ^{2}}{2h}\\
&  \ \ \ \ -\frac{1}{h}\sin^{2}\left(  \sigma\right)  \left(  Du(x_{h}%
)-Du(x_{0})\right)  \cdot\left(  x_{0}+h\vec{V}\right)  -\frac{1}{h}\sin
^{2}\left(  \sigma\right)  Du(x_{0})\cdot\left(  \left(  x_{0}+h\vec
{V}\right)  -x_{0}\right) \\
&  =\frac{u(x_{0}+h\vec{V})-u(x_{0})}{h}-\sin^{2}\left(  \sigma\right)
Du(x_{0})\cdot\vec{V}\\
&  \ \ \ \ +\cos\left(  \sigma\right)  \frac{\left[  \sin\left(
\sigma\right)  \left(  Du(x_{0}+h\vec{V})-Du(x_{0})\right)  \right]  \left[
Du(x_{0}+h\vec{V})+Du(x_{0})\right]  }{2h}\\
&  \ \ \ \ -\sin\left(  \sigma\right)  \cos\left(  \sigma\right)  \left(
x_{0}\cdot\vec{V}+\frac{h}{2}\left\vert \vec{V}\right\vert ^{2}\right) \\
&  \ \ \ \ -\frac{1}{h}\sin\left(  \sigma\right)  \left[  \sin\left(
\sigma\right)  \left(  Du(x_{0}+h\vec{V})-Du(x_{0})\right)  \right]
\cdot\left(  x_{0}+h\vec{V}\right)  .
\end{align*}
Rewriting (\ref{differenceh}) as
\begin{equation}
\sin\left(  \sigma\right)  \left[  Du(x_{h})-Du(x_{0})\right]  =h\bar{e}%
_{j}-\cos\left(  \sigma\right)  h\vec{V}%
\end{equation}
we see
\begin{align*}
\frac{\bar{u}\left(  \bar{x}_{h}\right)  -\bar{u}\left(  \bar{x}_{0}\right)
}{h}  &  =\frac{u(x_{0}+h\vec{V})-u(x_{0})}{h}-\sin^{2}\left(  \sigma\right)
Du(x_{0})\cdot\vec{V}\\
&  +\cos\left(  \sigma\right)  \frac{\left[  h\bar{e}_{j}-\cos\left(
\sigma\right)  h\vec{V}\right]  \left[  Du(x_{0}+h\vec{V})+Du(x_{0})\right]
}{2h}\\
&  -\sin\left(  \sigma\right)  \cos\left(  \sigma\right)  \left(  x_{0}%
\cdot\vec{V}+\frac{h}{2}\left\vert \vec{V}\right\vert ^{2}\right)  -\frac
{1}{h}\sin\left(  \sigma\right)  \left[  h\bar{e}_{j}-\cos\left(
\sigma\right)  h\vec{V}\right]  \cdot\left(  x_{0}+h\vec{V}\right) \\
&  =\frac{u(x_{0}+h\vec{V})-u(x_{0})}{h}-\sin^{2}\left(  \sigma\right)
Du(x_{0})\cdot\vec{V}\\
&  +\cos\left(  \sigma\right)  \frac{1}{2}\left[  \bar{e}_{j}-\cos\left(
\sigma\right)  \vec{V}\right]  \left[  2Du(x_{0})+\frac{h\bar{e}_{j}%
-\cos\left(  \sigma\right)  h\vec{V}}{\sin\left(  \sigma\right)  }\right] \\
&  -\sin\left(  \sigma\right)  \cos\left(  \sigma\right)  \left(  x_{0}%
\cdot\vec{V}+\frac{h}{2}\left\vert \vec{V}\right\vert ^{2}\right)
-\sin\left(  \sigma\right)  \left[  \bar{e}_{j}-\cos\left(  \sigma\right)
\vec{V}\right]  \cdot\left(  x_{0}+h\vec{V}\right) \\
&  =\frac{u(x_{0}+h\vec{V})-u(x_{0})}{h}-\sin^{2}\left(  \sigma\right)
Du(x_{0})\cdot\vec{V}\\
&  +\cos\left(  \sigma\right)  \left[  \bar{e}_{j}-\cos\left(  \sigma\right)
\vec{V}\right]  \cdot Du(x_{0})+\frac{h}{2}\frac{\cos\left(  \sigma\right)
}{\sin\left(  \sigma\right)  }\left\vert \bar{e}_{j}-\cos\left(
\sigma\right)  \vec{V}\right\vert ^{2}\\
&  -\sin\left(  \sigma\right)  \cos\left(  \sigma\right)  x_{0}\cdot\vec
{V}-\sin\left(  \sigma\right)  \cos\left(  \sigma\right)  \frac{h}%
{2}\left\vert \vec{V}\right\vert ^{2}-\sin\left(  \sigma\right)  \bar{e}%
_{j}\cdot x_{0}-h\sin\left(  \sigma\right)  \bar{e}_{j}\cdot\vec{V}\\
&  +\sin\left(  \sigma\right)  \cos\left(  \sigma\right)  x_{0}\cdot\vec
{V}+h\sin\left(  \sigma\right)  \cos\left(  \sigma\right)  \left\vert \vec
{V}\right\vert ^{2}\\
&  =\frac{u(x_{0}+h\vec{V})-u(x_{0})}{h}-\sin^{2}\left(  \sigma\right)
Du(x_{0})\cdot\vec{V}\\
&  +\cos\left(  \sigma\right)  \bar{e}_{j}\cdot Du(x_{0})-\cos^{2}\left(
\sigma\right)  Du(x_{0})\cdot\vec{V}-\sin\left(  \sigma\right)  \bar{e}%
_{j}\cdot x_{0}\\
&  +h\left[
\begin{array}
[c]{c}%
\frac{\cos\left(  \sigma\right)  }{\sin\left(  \sigma\right)  }\frac{1}%
{2}\left\vert \bar{e}_{j}-\cos\left(  \sigma\right)  \vec{V}\right\vert
^{2}-\sin\left(  \sigma\right)  \cos\left(  \sigma\right)  \frac{1}%
{2}\left\vert \vec{V}\right\vert ^{2}\\
-\sin\left(  \sigma\right)  \bar{e}_{j}\cdot\vec{V}+\sin\left(  \sigma\right)
\cos\left(  \sigma\right)  \left\vert \vec{V}\right\vert ^{2}%
\end{array}
\right] \\
&  =Du(x^{\ast})\cdot V-Du(x_{0})\cdot\vec{V}+\cos\left(  \sigma\right)
\bar{e}_{j}\cdot Du(x)-\sin\left(  \sigma\right)  \bar{e}_{j}\cdot x_{0}\\
&  +h\left[
\begin{array}
[c]{c}%
\frac{\cos\left(  \sigma\right)  }{\sin\left(  \sigma\right)  }\frac{1}%
{2}\left\vert \bar{e}_{j}-\cos\left(  \sigma\right)  \vec{V}\right\vert
^{2}-\sin\left(  \sigma\right)  \cos\left(  \sigma\right)  \frac{1}%
{2}\left\vert \vec{V}\right\vert ^{2}\\
-\sin\left(  \sigma\right)  \bar{e}_{j}\cdot\vec{V}+\sin\left(  \sigma\right)
\cos\left(  \sigma\right)  \left\vert \vec{V}\right\vert ^{2}%
\end{array}
\right]
\end{align*}
where $x^{\ast}$ is some value between $x_{0}+h\vec{V}$ and $x_{0}$ obtained
by the mean value theorem. Now we may take a limit with $h$ vanishing. Because
$\vec{V}$ (which a priori can point in many directions) is bounded, the
$h$-term vanishes in the limit. Because $Du$ is continuous, and $x(\bar{x})$
is Lipschitz, we also have that
\[
\lim_{h\rightarrow0}\left\vert \left(  Du(x^{\ast})-Du(x_{0})\right)  \cdot
V\right\vert \leq\lim_{h\rightarrow0}\sup\left\vert Du(x^{\ast})-Du(x_{0}%
)\right\vert \left\vert V\right\vert =0.
\]

We are left with
\begin{equation}
\lim_{h\rightarrow0}\frac{\bar{u}(\bar{x}_{0}+h\bar{e}_{j})-\bar{u}(\bar
{x}_{0})}{h}=\cos\left(  \sigma\right)  u_{j}(x_{0})-\sin\left(
\sigma\right)  x_{0}^{j}. \label{star}%
\end{equation}
This is precisely the $\bar{y}$-component of the image of the rotation
(\ref{xbarybar}). It follows that the gradient graph of $\bar{u}$ exists
everywhere and is isometric to the gradient graph of $u.$
\end{proof}

\begin{cor}
\label{conjrots}\bigskip An analogous result holds when $u$ is semi-concave,
and $\sigma$ is negative. The rotations through $\sigma$ and $-\sigma$ are
inverse operations where they are defined, up to an additive constant in the
potential function.
\end{cor}

\begin{proof}
While we could claim a proof that is formally the same as the proof of
Proposition \ref{weakrotation}, we offer an alternative argument based on the
fact that, whenever $u$ is semi-concave, $-u$ must be semi-convex. \ Starting
with a semi-convex $-u,$ we may rotate the graph $\Gamma_{-u}$ by a downward
rotation through $-\sigma,$ applying Proposition \ref{weakrotation}, and then
take the complex conjugate of the result in $\mathbb{C}^{n}$. This follows
from the fact that, as operators on $\mathbb{C}^{n}$ ($\mathbb{R}$-linear on
$\mathbb{R}^{2n})$ for any diagonal unitary matrix $U$ we have
\[
c\circ U\circ c=U^{-1}=U^{\ast}%
\]
where $c$ is the $\mathbb{R}$-linear complex conjugation map on $\mathbb{R}%
^{2n},$ that is
\[
c(x+\sqrt{-1}y)=x-\sqrt{-1}y.
\]
In particular, taking $-\overline{(-u)}$ via rotation of $-u$ (not complex
conjugation), we obtain the potential $\bar{u}$ for the graph rotated through
a negative angle $-\sigma.$
\end{proof}

The following technical result is useful when we approximate $u$ while keeping
$K$-convexity.

\begin{lem}
\label{convexmoll}\ Let $u^{\varepsilon\text{ }}$ be a standard mollification
of $u.$ \ If $u$ is $K$-convex on $\Omega$, then so is $u^{\varepsilon}$ on
\begin{equation}
\Omega^{\varepsilon}=\left\{  x:d(x,\partial\Omega)>\varepsilon\right\}  .
\label{Omegaeps}%
\end{equation}

\end{lem}

\begin{proof}
Consider a mollifier $\phi$ that is radial, supported in $B_{\varepsilon
}\left(  0\right)  $ and has unit integral. Given a point $x\in$
$\Omega^{\varepsilon}$,
\begin{align*}
u^{\varepsilon}(x)  &  =\int_{\Omega}\phi(x-y)u(y)dy\\
&  =\int_{B_{\varepsilon}(x)}\phi(x-y)u(y)dy\\
&  =\int_{B_{\varepsilon}(0)}\phi(z)u(x+z)dz
\end{align*}
so we have
\[
Du^{\varepsilon}(x)=\int_{B_{\varepsilon}(0)}\phi(z)Du(x+z)dz
\]
Now consider, for $x_{1},x_{0}\in$ $\Omega^{\varepsilon},$ the expression
\begin{align*}
&  \langle Du^{\varepsilon}(x_{1})-Du^{\varepsilon}(x_{0}),x_{1}-x_{0}%
\rangle\\
&  =\left\langle \int_{B_{\varepsilon}(0)}\phi(z)Du(x_{1}+z)dz-\int%
_{B_{\varepsilon}(0)}\phi(z)Du(x_{0}+z)dz,x_{1}-x_{0}\right\rangle \\
&  =\left\langle \int_{B_{\varepsilon}(0)}\phi(z)\left(  Du(x_{1}%
+z)-Du(x_{0}+z)\right)  ,x_{1}-x_{0}\right\rangle dz\\
&  =\int_{B_{\varepsilon}(0)}\phi(z)\left\langle Du(x_{1}+z)-Du(x_{0}%
+z),(x_{1}+z)-\left(  x_{0}+z\right)  \right\rangle dz\\
&  \geq\int_{B_{\varepsilon}(0)}\phi(z)K\left\vert x_{1}-x_{0}\right\vert
^{2}dz\\
&  =K\left\vert x_{1}-x_{0}\right\vert ^{2}.
\end{align*}

\end{proof}

\begin{prop}
\label{convexityestimates}Suppose that $u$ is $\tan(\kappa)$-convex and
$C^{1}$ and $\bar{u}$ is obtained as in Proposition \ref{weakrotation}. If
$\kappa,\sigma,\kappa-\sigma\in\left(  -\pi/2,\pi/2 \right)  $, then $\bar{u}$
is $\tan(\kappa-\sigma)$-convex.
\end{prop}

\begin{proof}
We define the following functions
\begin{align*}
\bar{x}_{\varepsilon}  &  =\cos(\sigma)x+\sin\left(  \sigma\right)
Du^{\varepsilon}(x)\\
\bar{y}_{\varepsilon}  &  =-\sin\left(  \sigma\right)  x+\cos(\sigma
)Du^{\varepsilon}(x).
\end{align*}
Note that, as before, the set
\[
\bar{\Gamma}_{\varepsilon}=\left\{  \left(  \bar{x}_{\varepsilon}(x),\bar
{y}_{\varepsilon}(x)\right)  :x\in\Omega\right\}
\]
is the rotation of the gradient graph of $u^{\varepsilon}$ through angle
$\sigma.$ \ \ (To be clear, we are not taking the gradient graph of the
mollified rotated function, rather we are rotating the gradient graph of the
mollified function.)\ 

Now $Du$ is continuous, so the mollified derivatives $Du^{\varepsilon}$ will
converge locally uniformly to $Du$ as $\varepsilon$ $\rightarrow0$ (cf.
\cite[Appendix C, Theorem 6]{Evans}). \ It follows that the functions $\bar
{x}_{\varepsilon}$ and $\bar{y}_{\varepsilon}$ will also converge locally
uniformly, to $\bar{x}$ and $\bar{y}$ respectively, as functions of $x,$
where
\begin{align*}
\bar{x}  &  =\cos(\sigma)x+\sin\left(  \sigma\right)  Du(x)\\
\bar{y}  &  =-\sin\left(  \sigma\right)  x+\cos(\sigma)Du(x).
\end{align*}
We have seen in Proposition \ref{weakrotation} that
\[
\bar{\Gamma}=\left\{  \left(  \bar{x}(x),\bar{y}(x)\right)  :x\in
\Omega\right\}
\]
is precisely the gradient graph of the function $\bar{u}$ over $\bar{\Omega}.$
The semi-convexity condition (\ref{Kconvex}) on $\bar{u}$ that we are trying
to show is
\[
\langle\bar{y}(x_{1})-\bar{y}(x_{0}),\bar{x}(x_{1})-\bar{x}(x_{0})\rangle
\geq\tan(\kappa-\sigma)\left\vert \bar{x}(x_{1})-\bar{x}(x_{0})\right\vert
^{2}.
\]
We claim that%
\begin{equation}
\langle\bar{y}_{\varepsilon}(x_{1})-\bar{y}(x_{0}),\bar{x}_{\varepsilon}%
(x_{1})-\bar{x}(x_{0})\rangle\geq\tan(\kappa-\sigma)\left\vert \bar
{x}_{\varepsilon}(x_{1})-\bar{x}(x_{0})\right\vert ^{2} \label{barconvex}%
\end{equation}
for all $\varepsilon>0$. The local uniform convergence of $\bar{x}%
_{\varepsilon}$ and $\bar{y}_{\varepsilon}$ will then give us the result. \ To
show (\ref{barconvex}), we start by computing the Jacobian of the map $\bar
{x}_{\varepsilon}:$

Since $u^{\varepsilon}$ is smooth
\[
\frac{d\bar{x}_{\varepsilon}}{dx}=\cos(\sigma)I_{n}+\sin\left(  \sigma\right)
D^{2}u^{\varepsilon}(x).
\]
By assumption, $u$ is $\tan\left(  \kappa\right)  $-convex, and hence so is
$u^{\varepsilon}$, by Lemma \ref{convexmoll}, at least on $\Omega
^{\varepsilon}$ (recall (\ref{Omegaeps})). It follows that
\[
D^{2}u^{\varepsilon}(x)\geq\tan\left(  \kappa\right)  I_{n}.
\]
So%
\begin{align*}
\frac{d\bar{x}_{\varepsilon}}{dx}  &  \geq\cos(\sigma)I_{n}+\sin\left(
\sigma\right)  \tan\left(  \kappa\right)  I_{n}\\
&  =\frac{\cos(\sigma-\kappa)}{\cos\left(  \kappa\right)  }I_{n}>0
\end{align*}
since $\kappa$ and $\sigma-k\in(-\pi/2,\pi/2)$. The coordinate change is
invertible and the Jacobian can be computed
\[
\frac{dx}{d\bar{x}_{\varepsilon}}=\left(  \cos(\sigma)I_{n}+\sin\left(
\sigma\right)  D^{2}u^{\varepsilon}(x)\right)  ^{-1}.
\]
Next
\[
D\bar{y}_{\varepsilon}=\left(  -\sin\left(  \sigma\right)  I_{n}+\cos
(\sigma)D^{2}u^{\varepsilon}(x)\right)  .
\]
Now each $\bar{\Gamma}_{\varepsilon}$ is the gradient graph of a function
$\bar{u}_{\varepsilon}\left(  \bar{x}_{\varepsilon}\right)  $ on the region
$\bar{x}_{\varepsilon}\left(  \Omega\right)  $. In order to compute the
Hessian of $\bar{u}_{\varepsilon}$ in terms of $\bar{x}_{\varepsilon},$ we
compute%
\begin{align*}
D_{\bar{x}_{\varepsilon}}^{2}\bar{u}_{\varepsilon}  &  =D_{x}\bar
{y}_{\varepsilon}\cdot\frac{dx}{d\bar{x}_{\varepsilon}}=D_{\bar{x}%
_{\varepsilon}}\bar{y}_{\varepsilon}\\
&  =\left(  -\sin\left(  \sigma\right)  I_{n}+\cos(\sigma)D^{2}u^{\varepsilon
}(x)\right)  \left(  \cos(\sigma)I_{n}+\sin\left(  \sigma\right)
D^{2}u^{\varepsilon}(x)\right)  ^{-1}.
\end{align*}
At any point, we may diagonalize the expression for $D_{\bar{x}^{\varepsilon}%
}^{2}\bar{u}_{\varepsilon}(\bar{x})$ by diagonalizing $D^{2}u^{\varepsilon
}(x(\bar{x}))$: \
\[
D_{\bar{x}^{\varepsilon}}^{2}\bar{u}_{\varepsilon}=\left(
\begin{array}
[c]{ccc}%
\frac{-\sin\left(  \sigma\right)  +\cos\left(  \sigma\right)  \lambda_{1}%
}{\cos\left(  \sigma\right)  +\sin\left(  \sigma\right)  \lambda_{1}} & 0 &
0\\
0 & ... & 0\\
0 & 0 & \frac{-\sin\left(  \sigma\right)  +\cos\left(  \sigma\right)
\lambda_{n}}{\cos\left(  \sigma\right)  +\sin\left(  \sigma\right)
\lambda_{n}}%
\end{array}
\right)  =\left(
\begin{array}
[c]{ccc}%
\bar{\lambda}_{1} & 0 & 0\\
0 & ... & 0\\
0 & 0 & \bar{\lambda}_{n}%
\end{array}
\right)  .
\]
Now
\[
\bar{\lambda}_{j}=\frac{-\sin\left(  \sigma\right)  +\cos\left(
\sigma\right)  \lambda_{j}}{\cos\left(  \sigma\right)  +\sin\left(
\sigma\right)  \lambda_{j}}=\frac{-\frac{\sin\left(  \sigma\right)  }%
{\cos\left(  \sigma\right)  }+\lambda_{j}}{1+\frac{\sin\left(  \sigma\right)
}{\cos\left(  \sigma\right)  }\lambda_{j}}=\tan(-\sigma+\arctan(\lambda
_{j})).
\]
\ Because
\[
\arctan(\lambda_{j})\geq\kappa
\]
we conclude that
\[
\bar{\lambda}_{j}\geq\tan(-\sigma+\kappa)
\]
and $D_{\bar{x}^{\varepsilon}}^{2}\bar{u}_{\varepsilon}$ is $\tan
(-\sigma+\kappa)$-convex, that is
\begin{equation}
\langle D_{\bar{x}_{\varepsilon}}\bar{u}_{\varepsilon}(x_{1})-D_{\bar
{x}_{\varepsilon}}\bar{u}_{\varepsilon}(x_{0}),\bar{x}_{\varepsilon}%
(x_{1})-\bar{x}_{\varepsilon}(x_{0})\rangle\geq\tan(-\sigma+\kappa)\left\vert
\bar{x}_{\varepsilon}(x_{1})-\bar{x}_{\varepsilon}(x_{0})\right\vert ^{2}%
\end{equation}
or
\begin{equation}
\langle\bar{y}_{\varepsilon}(x_{1})-\bar{y}_{\varepsilon}(x_{0}),\bar
{x}_{\varepsilon}(x_{1})-\bar{x}_{\varepsilon}(x_{0})\rangle\geq\tan
(-\sigma+\kappa)\left\vert \bar{x}_{\varepsilon}(x_{1})-\bar{x}_{\varepsilon
}(x_{0})\right\vert ^{2}%
\end{equation}
provided that $x_{1}$ and $x_{0}$ are at least $\varepsilon$ away from the
boundary of $\Omega$. By the local uniform convergence, we conclude that
\begin{equation}
\left\langle \bar{y}(x_{1})-\bar{y}(x_{0}),\bar{x}(x_{1})-\bar{x}%
(x_{1})\right\rangle \geq\tan(-\sigma+\kappa)\left\vert \bar{x}(x_{1})-\bar
{x}(x_{1})\right\vert ^{2}%
\end{equation}
that is, $\bar{u}$ is $\tan(\kappa-\sigma)$-convex.
\end{proof}

The following is an observation on how semi-convexity can lead to bounded
geometry, even when the potential is not twice differentiable.

\begin{cor}
Suppose that $u\in C^{1}$ and is semi-convex. Then the gradient graph of $u$
is isometric to the gradient graph of a $C^{1,1}$ function.
\end{cor}

\begin{proof}
Choose $\sigma$ $\in(0,\pi/2)$ and $\varepsilon>0$ for which
(\ref{semiconvexity}) is satisfied. Now to control the $C^{1,1}$ norm of
$\bar{u}$ we note that
\begin{align*}
\left\Vert \bar{u}\right\Vert _{C^{1,1}\left(  \bar{\Omega}\right)  }  &
=\sup_{\bar{x}_{0},\bar{x}_{1}\in\bar{\Omega}}\frac{\left\vert D\bar{u}%
(\bar{x}_{1})-D\bar{u}(\bar{x}_{0})\right\vert }{\left\vert \bar{x}_{1}%
-\bar{x}_{0}\right\vert }\\
&  =\sup_{x_{0},x_{1}\in\Omega}\frac{\left\vert \bar{y}(x_{1})-\bar{y}%
(x_{0})\right\vert }{\left\vert \bar{x}(x_{1})-\bar{x}(x_{0})\right\vert }.
\end{align*}

So for any pair $x_{0},x_{1}\in\Omega$
\begin{align*}
\frac{\left\vert \bar{y}(x_{1})-\bar{y}(x_{0})\right\vert }{\left\vert \bar
{x}(x_{1})-\bar{x}(x_{0})\right\vert }  &  =\frac{\left\vert \cos\left(
\sigma\right)  Du(x_{1})-\sin\left(  \sigma\right)  x_{1}-\cos\left(
\sigma\right)  Du(x_{0})-\sin\left(  \sigma\right)  x_{0}\right\vert
}{\left\vert \cos\left(  \sigma\right)  x_{1}+\sin\left(  \sigma\right)
Du(x_{1})-\cos\left(  \sigma\right)  x_{0}+\sin\left(  \sigma\right)
Du(x_{0})\right\vert }\\
&  =\frac{\left\vert \cos\left(  \sigma\right)  \left(  Du(x_{1}%
)-Du(x_{0})\right)  -\sin\left(  \sigma\right)  \left(  x_{1}-x_{0}\right)
\right\vert }{\left\vert \cos\left(  \sigma\right)  \left(  x_{1}%
-x_{0}\right)  +\sin\left(  \sigma\right)  \left(  Du(x_{1})-Du(x_{0})\right)
\right\vert }.
\end{align*}
To show this is bounded, we explore two cases. \ Let $A=2\cot(\sigma)>0$.
\ The first case is when
\begin{equation}
\left\vert Du(x_{1})-Du(x_{0})\right\vert \leq A\left\vert x_{1}%
-x_{0}\right\vert . \label{case1}%
\end{equation}
Recall $\sigma\in(0,\pi/2)$, we have
\[
\frac{\left\vert \cos\left(  \sigma\right)  \left(  Du(x_{1})-Du(x_{0}%
)\right)  -\sin\left(  \sigma\right)  \left(  x_{1}-x_{0}\right)  \right\vert
}{\left\vert \cos\left(  \sigma\right)  \left(  x_{1}-x_{0}\right)
+\sin\left(  \sigma\right)  \left(  Du(x_{1})-Du(x_{0})\right)  \right\vert
}\leq\frac{\left\vert \cos\left(  \sigma\right)  A\left\vert x_{1}%
-x_{0}\right\vert +\sin\left(  \sigma\right)  \left\vert x_{1}-x_{0}%
\right\vert \right\vert }{\left\vert \cos\left(  \sigma\right)  \left(
x_{1}-x_{0}\right)  +\sin\left(  \sigma\right)  \left(  Du(x_{1}%
)-Du(x_{0})\right)  \right\vert }%
\]
and
\begin{align*}
&  \left\langle \cos\left(  \sigma\right)  \left(  x_{1}-x_{0}\right)
+\sin\left(  \sigma\right)  \left(  Du(x_{1})-Du(x_{0})\right)  ,\frac
{x_{1}-x_{0}}{\left\vert x_{1}-x_{0}\right\vert }\right\rangle \\
&  =\cos\left(  \sigma\right)  \left\vert x_{1}-x_{0}\right\vert +\left\langle
\sin\left(  \sigma\right)  \left(  Du(x_{1})-Du(x_{0})\right)  ,\frac
{x_{1}-x_{0}}{\left\vert x_{1}-x_{0}\right\vert }\right\rangle \\
&  \geq\cos\left(  \sigma\right)  \left\vert x_{1}-x_{0}\right\vert
+\sin\left(  \sigma\right)  \left\vert x_{1}-x_{0}\right\vert \left(
-\cot(\sigma)+\varepsilon\right) \\
&  =\sin\left(  \sigma\right)  \left\vert x_{1}-x_{0}\right\vert \varepsilon
\end{align*}
where we used (\ref{econvexity}) in the second line. Thus (\ref{case1}) leads
to%
\[
\frac{\left\vert \bar{y}(x_{1})-\bar{y}(x_{0})\right\vert }{\left\vert \bar
{x}(x_{1})-\bar{x}(x_{0})\right\vert }\leq\left\vert \frac{\cos\left(
\sigma\right)  A+\sin\left(  \sigma\right)  }{\sin\left(  \sigma\right)
\varepsilon}\right\vert =\frac{\cos^{2}\left(  \sigma\right)  +1}{\sin
^{2}\left(  \sigma\right)  }\frac{1}{\varepsilon}.
\]
The next case is when
\begin{equation}
\left\vert Du(x_{1})-Du(x_{0})\right\vert \geq A\left\vert x_{1}%
-x_{0}\right\vert . \label{case2}%
\end{equation}
Then by the triangle inequality and \eqref{case2}
\begin{align*}
\left\vert \cos\left(  \sigma\right)  \left(  x_{1}-x_{0}\right)  +\sin\left(
\sigma\right)  \left(  Du(x_{1})-Du(x_{0})\right)  \right\vert  &  \geq
\sin(\sigma)|Du(x_{1})-Du(x_{0})|-\cos(\sigma)|x_{1}-x_{0}|\\
&  \geq\left(  \sin\left(  \sigma\right)  -\frac{\cos(\sigma)}{A}\right)
\left\vert Du(x_{1})-Du(x_{0})\right\vert \\
&  =\frac{1}{2}\sin\left(  \sigma\right)  \left\vert Du(x_{1})-Du(x_{0}%
)\right\vert
\end{align*}
and
\begin{align*}
\frac{\left\vert \cos\left(  \sigma\right)  \left(  Du(x_{1})-Du(x_{0}%
)\right)  -\sin\left(  \sigma\right)  \left(  x_{1}-x_{0}\right)  \right\vert
}{\left\vert \cos\left(  \sigma\right)  \left(  x_{1}-x_{0}\right)
+\sin\left(  \sigma\right)  \left(  Du(x_{1})-Du(x_{0})\right)  \right\vert }
&  \leq\frac{\cos\left(  \sigma\right)  \left(  Du(x_{1})-Du(x_{0})\right)
+\sin\left(  \sigma\right)  \frac{\left\vert Du(x_{1})-Du(x_{0})\right\vert
}{A}}{\frac{1}{2}\sin\left(  \sigma\right)  \left\vert Du(x_{1})-Du(x_{0}%
)\right\vert }\\
&  =\frac{\cos^{2}\left(  \sigma\right)  +1}{\sin\left(  \sigma\right)
\cos\left(  \sigma\right)  }.
\end{align*}
In either case, we have%
\[
\frac{\left\vert \bar{y}(x_{1})-\bar{y}(x_{0})\right\vert }{\left\vert \bar
{x}(x_{1})-\bar{x}(x_{0})\right\vert }\leq\max\left\{  \frac{\cos^{2}\left(
\sigma\right)  +1}{\sin^{2}\left(  \sigma\right)  }\frac{1}{\varepsilon}%
,\frac{\cos^{2}\left(  \sigma\right)  +1}{\sin\left(  \sigma\right)
\cos\left(  \sigma\right)  }\right\}  =C
\]
and $\bar{u}$ is $C^{1,1}.$\ 
\end{proof}

\bigskip The following corollary is immediate from the above by applying the
De Giorgi-Nash theorem.

\begin{cor}
Suppose that $u$ $\in C^{1\text{ }}$is a semi-convex weak solution to
(\ref{Hamstat}). Then the phase $\theta$ enjoys interior H\"{o}lder estimates
(with respect to the metric distances)\ on $\Gamma_{u}.$
\end{cor}

Finally, we show that smoothness and strong semi-concavity estimates on the
rotated potential can be used to conclude smoothness on $u$.

\begin{prop}
\label{smoothboth}\bigskip Suppose that $u$ and $\bar{u}$ are as in
Proposition \ref{weakrotation} and $\bar{u}\in$ $C^{2}\left(  \bar{\Omega
}\right)  $. Suppose also that for some constant $\epsilon>0$
\begin{equation}
D_{\bar{x}}^{2}\bar{u}\leq\left(  \frac{\cos\left(  \sigma\right)  }%
{\sin\left(  \sigma\right)  }-\epsilon\right)  I_{n}. \label{ubarbound}%
\end{equation}
Then for any integer $k>1$
\[
\left\Vert D^{k}u\right\Vert _{L^{\infty}(\Omega)}\leq C\left(  \sigma
,\epsilon,n\right)  \left(  \left\Vert D^{k}\bar{u}\right\Vert _{L^{\infty
}(\bar{\Omega})},\left\Vert D^{k-1}u\right\Vert _{L^{\infty}(\Omega)}\right)
.
\]

\end{prop}

\begin{proof}
The function $\bar{u}$ was obtained by a downward rotation of $\sigma$ from
$u,$ so $u$ may be obtained by the inverse rotation. In particular\ as
$\bar{u}\in$ $C^{2}\left(  \bar{\Omega}\right)  ,$ the change of variable
formulae hold on $\bar{\Omega}$:
\begin{align*}
x  &  = \cos(\sigma) \bar{x}-\sin(\sigma) D_{\bar{x}}\bar{u}(\bar{x})\\
y  &  = \sin(\sigma)\bar{x}+\cos(\sigma) D_{\bar{x}}\bar{u}(\bar{x}).
\end{align*}
Differentiating the first formula leads to
\[
\frac{dx}{d\bar{x}} =\cos\left(  \sigma\right)  I_{n}-\sin\left(
\sigma\right)  D^{2}_{\bar{x}}\bar{u}(\bar{x})\label{ijacobian}
\]
and noting that
\[
y =D_{x}u(x) =D_{x} u(x(\bar{x}))
\]
we have
\[
D_{x}u(\bar{x}) =\sin(\sigma)\bar{x}+\cos(\sigma)D_{\bar{x}}\bar{u}(\bar{x}).
\]
Now
\begin{align*}
D_{x}^{2}u  &  =D_{x}D_{x}u\\
&  =D_{x}\left(  \sin(\sigma)\bar{x}+\cos(\sigma)D_{\bar{x}}\bar{u}(\bar
{x})\right) \\
&  =\left(  \sin(\sigma)I_{n} +\cos(\sigma)D_{\bar{x}}^{2}\bar{u}(\bar
{x})\right)  \frac{d\bar{x}}{dx}.
\end{align*}
Noting (\ref{ubarbound}), we may invert (\ref{ijacobian}) and conclude%
\begin{align}
D_{x}^{2}u\left(  \bar{x}\right)   &  =\left(  \sin(\sigma)I+\cos
(\sigma)D_{\bar{x}}^{2}\bar{u}(\bar{x})\right)  \cdot\left(  \cos\left(
\sigma\right)  I_{n}-\sin\left(  \sigma\right)  D^{2}_{\bar{x}}\bar{u}(\bar
{x})\right)  ^{-1}\label{coh}\\
&  :=F_{\sigma}(D_{\bar{x}}^{2}\bar{u}(\bar{x}(x))).\nonumber
\end{align}

First, we will show that if $D_{\bar{x}}^{3}\bar{u}$ exists, then so will
$D_{x}^{3}u(x)$. To do this we differentiate (\ref{coh}) in $x$, obtaining
\begin{align*}
D_{x}D_{x}^{2}u(x)  &  =D_{x}F_{\sigma}(D_{\bar{x}}^{2}\bar{u}(\bar{x}(x)))\\
&  =\frac{dF_{\sigma}}{dD_{\bar{x}}^{2}\bar{u}}\cdot\frac{dD_{\bar{x}}^{2}%
\bar{u}}{d\bar{x}}\cdot\frac{d\bar{x}}{dx}.
\end{align*}
Combining (\ref{ubarbound}), the assumption that $D_{\bar{x}}^{3}\bar{u}$
exists, and the fact that all of these factors are well-defined and bounded,
we conclude that $D_{x}^{3}u$ exists and is controlled in terms of $D_{\bar
{x}}^{3}\bar{u}$.

Higher order estimates follow in the same way inductively.
\end{proof}

\section{Proof of Theorem \ref{T2}}

\begin{proof}
We are assuming that the function $\theta$ is a weak solution to a divergence
type equation (\ref{Hamstat}) on the set $\mathbb{B}_{1}(0)\backslash Q$.
Because the conditions (\ref{delta1}), (\ref{delta2}) and (\ref{delta3}) each
guarantee uniform ellipticity of the Laplace equation, we may immediately
apply Theorem \ref{Serrin} and conclude that $\theta$ is a weak solution over
the whole ball $\mathbb{B}_{1}(0)$.

Recall that
\[
F(D^{2}u)=F(\lambda_{1},\cdots,\lambda_{n})=\sum_{i=1}^{n}\arctan\lambda_{i}.
\]
To begin, we claim that if either of the conditions (\ref{delta1}) or
(\ref{delta2}) holds, then%
\[
F(D^{2}u)=\theta
\]
is a solution to a concave equation.

For the case $\theta\geq\delta+\frac{\pi}{2}(n-2)$, we recall that by
\cite[Lemma 2.1]{YY06} (see also \cite[section 8]{CNS3}) the level sets of
$f$, at any level $c$ with $|c|\geq\frac{\pi}{2}(n-2)$, are convex. \ We have
a uniform bound $\left\vert D^{2}u\right\vert \leq C_{0}$ wherever the Hessian
exists, so we may find a compact set ${\mathcal{K}}\subset S^{n\times n}$ such
that $F(M)>\frac{\pi}{2}(n-2)$ for any $M\in\mathcal{K}$, where $S^{n\times
n}$ is the space of symmetric $n\times n$ real matrices, such that
\begin{align*}
D^{2}u(\mathbb{B}_{1}(0))  &  \subset\mathcal{K}\text{ }\\
F(M)  &  >\frac{\delta}{2}+\frac{\pi}{2}(n-2)\text{ for all }M\in\mathcal{K}%
\end{align*}
We may smoothly modify $F$ on $\mathcal{K}$,
\[
\tilde{F}=f(F)
\]
so that $\tilde{F}$ is a uniformly concave function and has the same level
sets as $F$ on $\mathcal{K}$. (For a recent detailed proof of this fact, see
\cite[Lemma 2.2]{CPW} .) \ In this case
\[
\tilde{F}(D^{2}u)=\tilde{\theta}%
\]
for some smoothly modified $\tilde{\theta}$, constructed from $f$ such that
\[
\left\Vert \tilde{\theta}\right\Vert _{C^{\alpha}}\leq C\left\Vert
\theta\right\Vert _{C^{\alpha}}.
\]

For the second case, (\ref{delta2}), $u$ is uniformly convex, and the function
$F$ is clearly concave in the eigenvalues. So by taking $\tilde{F}=F$ (see
\cite[section 3]{CNS3}) we already have that%
\[
\tilde{F}(D^{2}u)=\theta
\]
for some concave $\tilde{F}$. Again, because $\left\vert D^{2}u\right\vert
\leq C_{0}$ where it exists, we can find a compact set $\mathcal{K}$ (still
using the same notation as above for simplicity) such that $D^{2}%
u(\mathbb{B}_{1}(0))\subset\mathcal{K}$ and $F$ is uniformly concave on
$\mathcal{K}$.

In either case, (\ref{delta1}) or (\ref{delta2}), we may extend $\tilde{F}$
beyond $\mathcal{K}$ to a global function $\bar{F}$ on $S^{n\times n}$ to
obtain a uniformly elliptic $\bar{F}$, satisfying $\bar{F}(M)=\tilde{F}(M)$
for $M\in\mathcal{K},$ $\bar{F}$ is uniformly elliptic, $\bar{F}$ is concave,
and $\bar{F}$ is continuous on $S^{n\times n}$ and still smooth on the
interior of $\mathcal{K}$. (For example, see \cite[Lemma 2.2]{Twisted}.)

Now we apply \cite[Theorem 8.1 and Remark 1 following, see also Remark 1 in
6.2]{CC}, which is Schauder theory for uniformly elliptic concave equations.
\ Note that \cite[p. 54 ]{CC} only requires the function $\bar{F}$ to be
concave and continuous. First note that by De Giorgi-Nash, when $u\in C^{1,1}$
the equation (\ref{Hamstat}) is uniformly elliptic, so the function
$\theta$ enjoys H\"{o}lder estimates. With the H\"{o}lder continuous function
$\theta$ determined, we can seek a viscosity solution of the boundary value
problem:
\begin{align*}
\bar{F}(D^{2}u^{\prime})  &  =\theta\text{ on }\mathbb{B}_{1}(0)\\
u^{\prime}  &  =u\text{ on }\partial\mathbb{B}_{1}(0).
\end{align*}
The viscosity solution exists by Perron's method, and is unique \cite[Theorem
4.1]{CIL}. \ Now our definition of weak solution is that $F(D^{2}u)=\theta$
almost everywhere, so we may apply \cite[Corollary 3]{LionsBMP} to conclude
that $u$ is also a solution to $\bar{F}(D^{2}u)=\theta$. \ Thus $u^{\prime}=u$
and all statements about viscosity solutions in \cite{CC} will apply to $u.$
Because the modification of $F$ was either smooth or away from a compact set
containing the image of $D^{2}u$, we still have
\[
\left\Vert \bar{F}(D^{2}u)\right\Vert _{C^{\alpha}\left(  \mathbb{B}%
_{4/5}(0)\right)  }\leq C_{1}%
\]
for some $C_{1}$ depending on the ellipticity constants, following from De
Giorgi-Nash, noting that $\left\Vert \theta\right\Vert _{L^{\infty}}\leq
n\pi/2.$ We conclude from \cite{CC} that
\[
\left\Vert D^{2}u\right\Vert _{C^{\alpha}\left(  \mathbb{B}_{3/4}(0)\right)
}\leq C_{2}%
\]
for $C_{2}$ depending on the ellipticity constants, $C_{1},$ and the
oscillation of $u.$ \ \ Now $\theta$ is a solution to a divergence type
equation with $C^{\alpha}$ coefficients, so we may apply \cite[Theorem
3.13]{HanLin} to conclude that
\[
\left\Vert \theta\right\Vert _{C^{1,\alpha}\left(  \mathbb{B}_{2/3}(0)\right)
}\leq C_{3}.
\]

Now for $e_{k\text{ }}$, consider the function
\[
\theta^{(h_{k})}(x)=\frac{\theta(x+he_{k\text{ }})-\theta(x)}{h}%
\]
defined on some interior region, for small $h>0.$ Because $\theta\in
C^{1,\alpha}\left(  \mathbb{B}_{2/3}(0)\right)  $ we have
\[
\left\Vert \theta^{(h_{k})}\right\Vert _{C^{\alpha}\left(  \mathbb{B}%
_{2/3-h}(0)\right)  }\leq C_{3}.
\]

Now
\begin{align*}
\theta^{(h_{k})}(x)  &  =\frac{1}{h}\int_{0}^{1}\frac{d}{dt}F(D^{2}%
u(x+he_{k\text{ }})t+(1-t)D^{2}u(x))dt\\
&  =\frac{1}{h}\int_{0}^{1}g^{ij}\left(  D^{2}u(x+he_{k\text{ }}%
)t+(1-t)D^{2}u(x)\right)  \left(  u(x+he_{k\text{ }})_{ij}-u_{ij}(x)\right)
dt\\
&  =\int_{0}^{1}g^{ij}\left(  D^{2}u(x+he_{k\text{ }})t+(1-t)D^{2}u(x)\right)
\left(  \frac{u(x+he_{k\text{ }})_{ij}-u_{ij}(x)}{h}\right)  dt\\
&  =G^{ij}u_{ij}^{(h_{k})}(x)\\
&  :=Lu^{(h_{k})}(x)
\end{align*}
for some uniformly elliptic $L=G^{ij}\partial_{i}\partial_{j}$ which is an
average of elliptic operators with $C^{\alpha}$ coefficients. Thus, each
$u^{(h_{k})}$ satisfies an uniformly elliptic equation of non-divergence type,
that is
\[
Lu^{(h_{k})}=\theta^{(h_{k})}\in C^{\alpha}\left(  \mathbb{B}_{2/3-h}%
(0)\right)
\]
with H\"{o}lder estimate uniform in $h$. Noting that each $u^{(h_{k})}\in
C^{2,\alpha}$ we may apply the non-divergence Schauder theory \cite[Theorem
6.6]{GT} to conclude a uniform $C^{2,\alpha}$ estimate as $h\rightarrow0$.
Thus, for each $k \in 1,...,n $ we have
\[
\left\Vert u_{k}\right\Vert _{C^{2,\alpha}\left(  \mathbb{B}_{1/2}(0)\right)
}\leq C_{4}%
\]
that is
\begin{align*}
u  &  \in C^{3,\alpha}\left(  \mathbb{B}_{1/2}(0)\right) \\
g  &  \in C^{1,\alpha}\left(  \mathbb{B}_{1/2}(0)\right)
\end{align*}
with estimates. \ 

Now from $\Delta_{g}\theta=0$ we get
\[
\sqrt{g}g^{ij}\theta_{ij}=-\partial_{i}\left(  \sqrt{g}g^{ij}\right)
\theta_{i}\in C^{\alpha}\left(  \mathbb{B}_{1/2}(0)\right)
\]
thus $\theta$ satisfies a non-divergence equation with H\"{o}lder continuous
right hand side. By Schauder theory \cite[Theorem 6.13]{GT}, $\theta$ must be
$C^{2,\alpha}.$ (More precisely, $\theta$ is the unique viscosity solution to
an equation which admits a $C^{2,\alpha}$ solution.) Iterating the previous
two steps, we may obtain all higher order estimates for any region further in
the interior.

Next we assume that (\ref{delta3}) holds. Suppose that a function $u$
satisfies (\ref{delta3}). Let
\[
\kappa=\arctan(1-\delta)<\frac{\pi}{4}.
\]
Condition (\ref{delta3}) gives us that $u$ is $-\tan\left(  \kappa\right)
$-convex. Perform a downward rotation of the graph of $u$ with $\sigma
=\frac{\pi}{4}.$ Proposition \ref{weakrotation} implies that the corresponding
coordinate change $\bar{x}(x)$ defined by (\ref{cc}) is bi-Lipschitz. It will
follow that any interior region of $\bar{\Omega}^{\varepsilon}$ (recall
\eqref{Omegaeps}) will be the homeomorphic image of an interior region
$\Omega^{\prime}$ with
\[
\Omega^{\varepsilon_{2}}\subset\Omega^{\prime}\subset\Omega^{\varepsilon_{1}}%
\]
with $\varepsilon_{1}/\varepsilon$ and $\varepsilon_{2}/\varepsilon$ bounded
above and away from $0$. It follows that interior estimates for $\bar{u}$ on
$\bar{\Omega}$ will correspond to interior estimates for $u$ on $\Omega.$

Now by Proposition \ref{convexityestimates}, $\bar{u}$ is $\beta_{0}$-convex
for
\[
\beta_{0}=\tan\left(  \arctan(\delta-1)-\frac{\pi}{4}\right)  =\frac{\delta
-2}{\delta}.
\]
Now letting $v=-u,$ we may also rotate \textit{upward} by $\sigma=\frac{\pi
}{4},$ to obtain a function $\bar{v}$ that is $\beta_{1}$-convex for
\[
\beta_{1}=\tan\left(  \arctan(\delta-1)+\frac{\pi}{4}\right)  =\frac{\delta
}{2-\delta}%
\]
by Proposition \ref{convexityestimates}. From the discussion in the proof of
Corollary \ref{conjrots}, we have that $\bar{v}=-\bar{u}$. In particular,
$-\bar{u}$ is $C^{1,1}$, uniformly convex, and clearly is also a weak solution
of (\ref{Hamstat}), as the quantity $\theta$ is odd in $D^{2}u$. We are then
back to the case (\ref{delta2}) , and may conclude interior estimates on the
derivatives of $-\bar{u}$ for any order, and hence also for derivatives of
$\bar{u}$. \ Now certainly (\ref{ubarbound}) holds for $\epsilon=1,$ so we may
apply Proposition \ref{smoothboth} and get interior derivative estimates on
$u$.
\end{proof}

\subsection{Proof of Theorem \ref{T3}}

\begin{proof}
Let $u$ be a $W^{2,n}\left(  \Omega\right)  $ solution to (\ref{vHamstat}).
Let $\Gamma_{u}=\left\{  (x,Du(x)):x\in\Omega\right\}  $. First note that the
Grassmannian geometry (in particular, the distance function) is invariant
under unitary actions on $\mathbb{C}^{n}$. Observe also that for small enough
$c_{0}\left(  n\right)  $, all Lagrangian planes within distance $c_{0}\left(
n\right)  $ from each other must be graphical over each other. Thus at any
point $p$ where $D^{2}u$ exists, the tangent space to $\Gamma$ is
well-defined, and we can locally take $\Gamma$ to be a graph over $T_{p}L.$ By
taking a unitary map sending $T_{p}\Gamma$ to $\mathbb{R}^{n}\times\{0\}$, we
may express the isometric image $\bar{\Gamma}$ locally as a gradient graph of
some function $\bar{u}$ over a region $\bar{\Omega}\subset\mathbb{R}^{n}$,
with $D^{2}\bar{u}(p)=0$. For Lagrangian tangent planes near $\mathbb{R}%
^{n}\times\left\{  0\right\}  ,$ the topology on the Lagrangian Grassmannian
is equivalent to the topology on Hessian space, so by choosing $c_{0}\left(
n\right)  $ small we have also guaranteed that
\[
\left\Vert u\right\Vert _{C^{1,1}(\Omega)}\leq c(n)<1
\]
where $c(n)$ is from Theorem \ref{T1}. Applying Theorem \ref{T1}, we may
conclude that $u$ is a weak solution to (\ref{Hamstat}). By Theorem \ref{T2},
$\bar{u}$ is smooth inside $\bar{\Omega}.$ So $\bar{\Gamma}$ is the gradient
graph of a smooth function over $\bar{\Omega}$, hence it is a smooth
submanifold of $%
\mathbb{R}
^{2n}.$
\end{proof}

Our result allows for the Hessian of the potential function $u$ to be just
continuous or even have mild discontinuities provided that $\Vert
u\Vert_{C^{1,1}}\leq c(n)$. The following result is obtained by Schoen and
Wolfson \cite[Proposition 4.6]{SWJDG}, for Lagrangian stationary surfaces
(when the potential functions are locally in $C^{2,\alpha}$) in general
K\"{a}hlerian ambient manifolds.

\begin{cor}
\label{C^2} Suppose that $u$ $\in C^{2}$ is a weak solution to (\ref{vHamstat}%
). Then $u$ is smooth.
\end{cor}

\begin{proof}
Let $\Gamma=\{\left(  x,Du(x)\right)  :x\in\Omega\}$. Near any point $x_{0}%
\in\Gamma$, we may write $\Gamma$ locally as as gradient graph of a function
$v$ over its tangent plane $T_{x_{0}}\Gamma$. \ Necessarily, this choice gives
us $D^{2}v(0)=0.$ \ Now $v$ is also stationary for compactly supported
variations near $x_{0}$, so $v$ must satisfy (\ref{vHamstat}) as well. Because
$D^{2}u\in C^{0}$, the tangent planes change continuously. It follows that
also $D^{2}v\in C^{0},$ and because we have chosen $D^{2}v(0)=0,$ we may find
a small neighborhood for which
\[
\left\Vert D^{2}v\right\Vert _{C^{0}}\leq c(n).
\]
Applying Theorem \ref{T3}, $v$ is smooth near $x$. It follows that $\Gamma$ is
smooth near $x$. \ Now because $D^{2}u$ was bounded, we may project the smooth
object $\Gamma$ back to the original coordinates $\Omega$, and the Jacobian
does not vanish. Thus we conclude that $u$ is a smooth function on $\Omega$.
\end{proof}

\bibliographystyle{amsalpha}
\bibliography{hamstat1}

\end{document}